\documentclass[12pt]{article}
\usepackage{amsfonts}
\usepackage{amsmath}
\usepackage{amssymb}
\usepackage{color}
\usepackage{latexsym}
\usepackage{mathrsfs} 
\usepackage{multirow}

\newcommand{\bc}{\begin{center}}
\newcommand{\ec}{\end{center}}
\newcommand{\ba}{\begin{array}}
\newcommand{\ea}{\end{array}}
\newcommand{\be}{\begin{eqnarray}}
\newcommand{\ee}{\end{eqnarray}}
\newcommand{\bel}{\begin{eqnarray}\label}
\newcommand{\eel}{\end{eqnarray}}
\newcommand{\bes}{\begin{eqnarray*}}
\newcommand{\ees}{\end{eqnarray*}}
\newcommand{\bi}{\begin{itemize}}
\newcommand{\ei}{\end{itemize}}
\newcommand{\bn}{\begin{enumerate}}
\newcommand{\en}{\end{enumerate}}
\newcommand{\pa}{\partial}
\newcommand{\f}{\frac}

\def\eps{\varepsilon}

\def\htheta{\hat{\theta}}
\def\btheta{\boldsymbol{\theta}}
\def\hbtheta{\hat{\btheta}}

\def\hlambda{\hat{\lambda}}
\def\bdelta{\boldsymbol{\delta}}

\def\hxi{\hat{\xi}}

\def\Gbar{\bar{G}}
\def\bI{\boldsymbol{I}}

\def\bu{\boldsymbol{u}}
\def\bv{\boldsymbol{v}}

\def\bx{\boldsymbol{x}}
\def\bX{\boldsymbol{X}}

\def\bZ{\boldsymbol{Z}}

\def\scrD{\mathscr{D}}

\def\scrR{\mathscr{R}}

\def\real{\mathbb{R}}

\def\argmin{\mathop{\rm arg\, min}}
\def\sgn{\mathrm{sgn}}
\def\Var{\mathrm{Var}}
\def\lam{{\lambda}}
\def\hlam{{\hlambda}}
\def\tlam{{\tilde \lambda}}

\newtheorem{thm}{Theorem}
\newtheorem{cor}{Corollary}
\newtheorem{lem}{Lemma}
\newtheorem{proposition}{Proposition}
\newtheorem{example}{Example}

\begin{document}


%
%
%
%
%
%
%
%

\centerline{\bf\large ADAPTIVE THRESHOLD ESTIMATION BY FDR}
\bigskip

\centerline{\sc Wenhua Jiang and Cun-Hui Zhang}
\medskip

\centerline{Soochow University and Rutgers University}
\bigskip

\begin{abstract}
This paper addresses the following simple question about sparsity.
For the estimation of an $n$-dimensional mean vector $\boldsymbol{\theta}$ in the Gaussian sequence model,
is it possible to find an adaptive optimal threshold estimator in a  full range of sparsity levels
where nonadaptive optimality can be achieved by threshold estimators?
We provide an explicit affirmative answer as follows.
Under the squared loss, adaptive minimaxity in strong and weak $\ell_p$ balls with $0\le p<2$ is
achieved by a class of smooth threshold estimators with the threshold level of the Benjamini-Hochberg
FDR rule or its a certain approximation, provided that the minimax risk is between
$n^{-\delta_n}$ and $\delta_n n$ for some $\delta_n\to 0$.
For $p=0$, this means adaptive minimaxity in $\ell_0$ balls when $1\le \|\boldsymbol{\theta}\|_0\ll n$.
The class of smooth threshold estimators includes the soft and firm threshold estimators
but not the hard threshold estimator.
The adaptive minimaxity in such a wide range is a delicate problem since the same is not true
for the FDR hard threshold estimator at certain threshold and nominal FDR levels.
The above adaptive minimaxity of the FDR smooth-threshold estimator is established
by proving a stronger notion of adaptive ratio optimality for the soft threshold estimator
in the sense that the risk for the FDR threshold level is uniformly within an infinitesimal fraction
of the risk for the optimal threshold level for each unknown vector, when the minimum risk of
nonadaptive soft threshold estimator is between $n^{-\delta_n}$ and $\delta_n n$.
It is an interesting consequence of this adaptive ratio optimality
that the FDR smooth-threshold estimator outperforms the sample mean
in the common mean model $\theta_i=\mu$ when $|\mu|<n^{-1/2}$.
\end{abstract}

%
%

\def\theequation{1.\arabic{equation}}
\setcounter{equation}{0}
\section{Introduction}\label{introduction}

Let $\bX=(X_1,\ldots,X_n)$ be a vector of independent variables with marginal
distributions $X_i\sim N(\theta_i,1)$, $i=1,\dots,n$.
Statistical inference about the mean vector $\btheta$, known as the Gaussian sequence problem,
has been considered as a canonical or motivating model in many important areas in statistics.
Examples include empirical Bayes, admissibility, nonparametric regression, variable
selection, multiple testing and so on.
It also carries immense practical relevance in statistical applications since observed
data are often understood, represented, or summarized approximately as a Gaussian vector.

An important Gaussian sequence problem, which illuminates our
understanding of more complicated models such as high-dimensional
linear regression and matrix estimation, is the estimation of a
sparse vector $\btheta$ under the squared-error loss. Donoho and
Johnstone \cite{DonohoJ94} made a fundamental contribution by
proving that when the sparsity of $\btheta$ is expressed as the
membership of $\btheta$ in a properly standardized small $\ell_p$
ball with $0\le p<2$, asymptotic minimaxity can be achieved by
threshold estimators but not linear estimators. Linear estimators do
not even achieve the optimal risk rate. However, since optimal
linear estimation in $\ell_2$ balls can be adaptively achieved by
the James-Stein estimator \cite{Stein56Berkeley, JamesS61}, it is a
natural question whether optimal threshold estimation can be
adaptively achieved with a data driven threshold level. The
universal threshold level $\sqrt{2\log n}$
\cite{DonohoJ94Biometrika}, equivalent to controlling the familywise
error rate in multiple testing, is suboptimal since it results in an
extra logarithmic risk factor for moderately small $\ell_p$ balls.

Adaptive threshold estimation has been considered by many, including
SURE \cite{DonohoJ95}, the generalized $C_p$ \cite{BirgeM01}, FDR
\cite{AbramovichBDJ06} and the parametric EB posterior median
(EBThresh) \cite{JohnstoneS04}. A general picture of the existing
analyses of these estimators is that adaptive exact minimax
threshold estimation is achieved in $\ell_p$ balls when the order of
the minimax risk is between $(\log n)^{1-p/2+\gamma}$ and
$n^{1-\kappa}$ with $\gamma>4.5$ and any $\kappa>0$
\cite{AbramovichBDJ06,WuZ13}, while adaptive rate minimaxity is
achieved when the minimax risk is of no smaller order than $O(1)$
\cite{BirgeM01}. We refer to Johnstone's book \cite{Johnstone11} for
a comprehensive discussion of the topic.

The state of the matter is not quite satisfactory in view of some
recent advances in empirical Bayes. In the framework of compound
estimation and empirical Bayes \cite{Robbins51, Robbins56}, our
problem can be considered as restricted or threshold empirical Bayes
\cite{Robbins83,Zhang03} since it aims to approximately achieve the
performance benchmark \bel{thresh-EB}
\min_{t(\cdot)\in \scrD}\int \int
\Big(t(x)-\theta\Big)^2\varphi(x-\theta)\,dx\,G_n(d\theta) 
\eel 
with a class of functions $\scrD$ restricted to threshold functions, where
$\varphi(z)$ is the $N(0,1)$ density and $G_n(u) =
n^{-1}\sum_{i=1}^n I\{\theta_i\le u\}$, the empirical Bayes nominal
prior, is the empirical distribution of the true deterministic
unknown mean vector $\btheta$. As we have mentioned, the James-Stein
estimator, which can be viewed as linear or parametric empirical
Bayes \cite{EfronM72,EfronM73,Morris83}, is not rate optimal for the
estimation of sparse $\btheta$. However, the general empirical Bayes
\cite{Robbins56,Robbins83}, which aims to approximately achieve the
benchmark (\ref{thresh-EB}) with general (unrestricted) $\scrD$,
enjoys optimality properties for sparse as well as dense signals
\cite{Zhang97, Zhang05, BrownG09, JiangZ09}. In fact, the general
maximum likelihood empirical Bayes (GMLEB) \cite{JiangZ09} is
guaranteed to possess the adaptive exact minimax and a stronger
adaptive ratio optimal properties when the order of the risk is
between $(\log n)^{5+3/p}$ and $n$ for sparse $\btheta$ and when
$n^{-1}\sum_{i=1}^n |\theta_i|^p\ll n^p/(\log n)^{4+9p/2}$ for dense
$\btheta$. Thus, as far as these first order optimality properties
are concerned, the advantage of threshold estimation is expected to
lie in a class of very sparse signals with risk of logarithmic
order, in view of the optimality of the GMLEB for both sparse and
dense signals. From this point of view, adaptive optimality
properties for widest range of sparse signals, especially for risks
at and below logarithmic rate, is highly desirable as a
theoretically nontrivial justification for the use of adaptive
threshold estimators when the signal is believed to be sparse.

The above discussion leads to the following interesting question.
When the signal is sparse, is the required lower bound for the risk for adaptive threshold estimation
of logarithmic order, of order 1, or even smaller?
Consider the case where minimax risk in $\ell_p$ balls is of smaller order than $n$
since threshold estimators do not asymptotically attain minimax risk anyway
in $\ell_p$ balls when the minimax risk is of order $n$.
In this case, the question can be phrased as whether adaptive optimal threshold estimation
can be achieved in $\ell_p$ balls when the minimax risk is above a logarithmic order, above
1 or smaller.

A main objective of this paper is to give an explicit affirmative
answer to the above question. We prove that with the threshold level
of the Benjamini-Hochberg rule \cite{BenjaminiH95} for controlling
the false discovery rate (FDR) or a suitable approximation of the
FDR rule, smooth threshold estimators between the soft and firm
threshold estimators uniformly achieve the adaptive exact minimaxity
in strong and weak $\ell_p$ balls when the minimax risk is between
$n^{-\delta_n}$ and $\delta_n n$ for any given $\delta_n\to 0$. For
$p=0$, this means adaptive minimaxity in $\ell_0$ balls when $1\le
\|\btheta\|_0\ll n$.

An interesting consequence of our result and that of
\cite{AbramovichBDJ06} is a proven advantage of smooth thresholding
against hard thresholding in adaptive estimation when the signal is
relatively weak and highly sparse. The importance of continuity was
advocated in \cite{FanL01} among others.

Another interesting question is whether it is possible for an adaptive threshold
estimator, designed to solve the nonparametric Gaussian sequence problem,
to also outperform a parametric estimator when the signal belongs to a parametric family.
Let ${\bar X} = n^{-1}\sum_{i=1}^nX_i$ be the sample mean.
In the common mean model $\theta_i=\mu\ \forall i\le n$, the $\ell_2$ risk of the sample mean
$\htheta_i = {\bar X}$ is $1$ for the estimation of $n$ elements of the vector $\btheta$.
Our results imply that in the common mean model, the FDR smooth threshold estimator
outperforms ${\bar X}$ when $|\mu|<(1-c)/\sqrt{n}$ for any fixed $c\in (0,1)$.
In fact, the FDR smooth threshold estimator is comparable to the optimal
soft thresholded ${\bar X}$ when $1/n^{1/2+\delta_n/2} \ll |\mu| \ll 1/n^{1/2}$.
See Example \ref{example-1} in Section 2.

These results are nontrivial in view of the following. Foster and
George \cite{ForsterG94} proved the existence of a risk inflation
factor of $2\log n$ between an optimal threshold estimator and an
oracle risk. Abramovich et al \cite{AbramovichBDJ06} proved that the
FDR hard threshold estimator is adaptive rate minimax when the order
of the minimax risk is between $(\log n)^{6-p/2}$ and $n^{1-\kappa}$
for all $\kappa>0$ but the same estimator is not adaptive minimax to
the constant factor when the nominal FDR level is higher than $1/2$.
Both results raised the possibility of a lower risk bound of
logarithmic order for adaptive estimation. The question below order
1 is even less clear. Birg\'e and Massart \cite{BirgeM01} raised the
possibility of a requirement of a risk of at least $O(1)$
($O(\eps^2)$ in their paper) for rate adaptive estimation of
$\btheta$. Moreover, when the minimax risk is of order $\log n$, the
nonadaptive asymptotic minimaxity of threshold estimators was only
proven recently \cite{Zhang12}.

Simultaneous adaptive threshold estimation for moderately and highly
sparse mean vectors is an analytically demanding problem.
Sophisticated machinery and clever arguments were deployed in
\cite{AbramovichBDJ06} to prove the theorem for minimax risk of
order between $(\log n)^{1-p/2+\gamma}$ and $n^{1-\kappa}$ with
$0\le p<2$, $\gamma\ge 5$ and $\kappa >0$. Their result was improved
upon recently to $\gamma>4.5$ in \cite{WuZ13}. We take a different
analytical approach by first proving a certain ratio optimality of
the FDR threshold level when the order of the minimum risk of soft
threshold estimation is between $n^{-\delta_n}$ and $\delta_n n$ at
the true unknown mean vector with $\delta_n\to 0$.

The ratio optimality asserts that the risk of the adaptive soft threshold estimator is uniformly within an
infinitesimal fraction of the minimum risk of soft threshold estimators over all threshold levels.
This directly guarantees the uniform optimality of the adaptive threshold level.
It is a stronger notion of optimality than adaptive minimaxity since it quarantees the performance
of the adaptive estimator at the true unknown $\btheta$ instead of the performance
in the worst case scenario in a class of the unknown $\btheta$.

The paper is organized as follows.
A class of FDR smooth threshold estimators is described in Section \ref{main-results},
along with statements of its adaptive ratio optimality and mininaxity properties.
Some preliminary analytical results are presented in Section \ref{analysis},
including some properties of the Bayes risk of smooth threshold estimators and its approximation,
comparison of the FDR rule and its population version, and applications of
the Gaussian isoperimetric inequalities to smooth  thresholding at random threshold levels.
Section \ref{sec-oracle} provides an oracle inequality and optimality properties
for a more general class of threshold rules than those in Section \ref{main-results}.
Section \ref{discussion} contains some discussion.
Mathematical proofs are provided in the Appendix.

\def\theequation{2.\arabic{equation}}
\setcounter{equation}{0}
\section{Main results}\label{main-results}

Let $P_{\btheta}$ denote probability measures under which
$X_i$ are independent statistics with
\bel{X_i}
X_i\sim N(\theta_i,1),\quad i=1,\ldots,n,
\eel
where $\btheta=(\theta_1,\ldots,\theta_n)$ is an unknown signal vector.
In the vector notation, it is convenient to state (\ref{X_i}) as $\bX = (X_1,\ldots,X_n) \sim N(\btheta,\bI_n)$
with $\bI_n$ being the identity matrix in $\real^n$.

Our problem is to estimate $\btheta$ under
the mean squared error
\bel{MSE}
E_{\btheta}\|\hbtheta-\btheta\|^2
=E_{\btheta}\sum_{i=1}^n(\htheta_i-\theta_i)^2\eel
for any estimator $\hbtheta=(\htheta_1,\ldots,\htheta_n)$.

Throughout the paper, boldface letters denote vectors and matrices,
for example, $\bX=(X_1,\ldots,X_n)$,
$\varphi(x)=e^{-x^2/2}/\sqrt{2\pi}$, $\Phi(t)=\int_{-\infty}^t\varphi(x)dx$ and $\Phi^{-1}(t)$
denote the standard normal density, distribution and quantile functions,
$\|\bv\|=\sqrt{\sum_iv_i^2}$, $\|\bv\|_\infty=\max_i|v_i|$ and $\|\bv\|_0=\{i:v_i\neq 0\}$
denote the $\ell_p$ norm for vectors $\bv$ with components $v_i$,
$x\vee y=\max(x,y)$, $x\wedge y=\min(x,y)$, $x_+=x\vee0$, $a_n\ll b_n$
means $a_n=o(b_n)$ and $a_n\approx b_n$ means
$\lim_{n\to\infty}a_n/b_n=1$. Univariate functions are applied to
vectors per component. Thus,
$\hbtheta=t(\bX)$ means $\htheta_i=t(X_i)$, $i=1,\ldots,n$.

\subsection{Adaptive threshold estimation by FDR}\label{subsec-FDR}
Given a sequence of null hypotheses $H_1,\ldots,H_n$, the false
discovery rate (FDR) of a multiple testing method is defined as \bes
\hbox{FDR} = E \frac{\#\{\hbox{\,falsely rejected
hypotheses\,}\}}{1\vee \#\{\hbox{\,rejected hypothesis\,}\}}. \ees
Benjamini and Hochberg \cite{BenjaminiH95} advocated the use of FDR
to measure Type-I errors in multiple testing and proposed the
following rule to control the FDR. Suppose independent test
statistics with known null distribution are observed for testing the
$n$ null hypotheses. Let $p_{(1)}\le p_{(2)}\le\cdots\le p_{(n)}$ be
the ordered $p$-values and $H_{(1)},\ldots,H_{(n)}$ the
corresponding hypotheses. The Benjamini-Hochberg rule controls the
FDR at the level $qn_0/n$ by rejecting hypotheses
$H_{(1)},\ldots,H_{({\hat k})}$, where ${\hat k} =\max\big\{i\colon
p_{(k)}\le q k /n\big\}$ and $n_0$ is the number of true hypotheses.
Since $n_0$ is unknown, the quantity $q$ is treated as the nominal
FDR level for the Benjamini-Hochberg rule.

Let $0<\alpha_2\le\alpha_1<1$. Define candidate threshold levels
\bel{xi_jk}
\xi_{j,k}= - \Phi^{-1}\Big(\frac{\alpha_jk}{2n}\Big),\quad j=1,2,
\eel
$k=1,\ldots,n$, and data-driven threshold levels
\bel{xi}
\hxi_1=\min\Big\{\xi_{1,k}\colon N(\xi_{1,k})\ge k\Big\},\quad
\hxi_2=\max\Big\{\xi_{2,k}\colon N(\xi_{2,k+1})<k+1\Big\},
\eel
where $N(t)$ is the number of observations above threshold $t$,
\bel{N(t)}
N(t)=\sum_{i=1}^nI\big\{|X_i|\ge t\big\}.
\eel

Let $H_i: \theta_i=0$ and considered two-sided test statistics $|X_i|$.
Since $X_i\sim N(\theta_i,0)$, the ordered $p$-values are $p_{(k)}=2\Phi(-|X_{(k)}|)$,
where $|X_{(1)}|\ge\cdots\ge|X_{(n)}|$ are the ordered values of $|X_i|$.
Since $N(\xi_{1,k})\ge k$ if and only if $p_{(k)}\le 2\Phi(-\xi_{1,k})=\alpha_1k/n$, $\hxi_1$ is the
threshold level of the Benjamini-Hochberg rule for controlling the FDR at the nominal level $\alpha_1$.

Likewise, $N(\xi_{2,k})\ge k$ if and only if $|X_{(k)}| > \xi_{2,k}$.
However, unlike $\hxi_1$, which corresponds to a step-up rule,
$\hxi_2$ is the threshold level of a step-down rule \cite{Holm79}
matched by the Benjamini-Hochberg rule.
Since $\xi_{1,k}\le \xi_{2,k}$, we have $\hxi_1 \le \hxi_2$ by (\ref{xi}).

Let $t_\lambda(x)$ denote a smooth threshold function indexed by its threshold level $\lam$.
We study optimality properties of the adaptive threshold estimator
\bel{hbtheta}
\hbtheta = t_{\hlambda}(\bX)= \big(t_{\hlambda}(X_1),\ldots,t_{\hlambda}(X_n)\big)
\eel
with a threshold level $\hlambda$ satisfying
\bel{hlam-simple}
\sqrt{1+\delta_{1,n}}\hxi_1 \le \hlambda \le \sqrt{1+\delta_{2,n}}\hxi_2,\quad 0\le\delta_{1,n}\le\delta_{2,n} \to 0.
\eel

The estimator (\ref{hbtheta}) is closely related to the $\ell_0$
penalized estimator \bel{PLSE} \hbtheta = \argmin_{\btheta}\left\{
\|\btheta - \bX\|_2^2 + \sum_{k=1}^{\|\btheta\|_0}\xi_k^2\right\}
\eel with nonincreasing threshold levels $\xi_k$. A local minimum of
(\ref{PLSE}), say $\hbtheta$ with $\|\hbtheta\|_0={\hat k}$, is a
hard threshold estimator at threshold level $\xi_{\hat k}$ such that
$X^2_{({\hat k})}\ge \xi_{\hat k}^2$ and $\xi_{{\hat k}+1}^2\ge
X^2_{({\hat k}+1)}$. When $\xi_k\in [\xi_{1,k}, \xi_{2,k}]$,
(\ref{xi}) implies $\hxi_1\le \xi_{\hat k} \le \hxi_2$ for such
local minima, so that (\ref{hlam-simple}) is satisfied with
$\delta_{1,n}=0$. In addition, (\ref{hbtheta}) allows the use of the
penalty level of (\ref{PLSE}) with data-dependent $\xi_k\in
[\xi_{1,k}, \xi_{2,k}]$. The $\ell_0$ penalized estimator
(\ref{PLSE}) was considered in
\cite{YangB98,ForsterS99,TibshiraniK99,GeorgeF00,BirgeM01,
AbramovichBDJ06,WuZ13} among many others.

\subsection{Adaptive ratio optimality of FDR rule in soft thresholding}\label{subsec-adaptive-ratio-opt}
Let $s_\lam(x)=\sgn(x)(|x|-\lam)_+$ be the soft threshold estimator.
Given a sequence of vector classes $\Theta^*_n\subset\real^n$, a threshold level $\hlambda$ is
adaptive ratio optimal for soft threshold estimation if
\bel{def-adaptive-ratio-opt}
\sup_{\btheta\in\Theta^*_n}\f{E_{\btheta}\|s_{\hlambda}(\bX)-\btheta\|^2}
{\inf_{\lambda\ge0}E_{\btheta}\|s_\lambda(\bX)-\btheta\|^2}\le1+o(1).
\eel
In words, the risk of $s_{\hlambda}$ is uniformly within an infinitesimal fraction
of the risk of $s_\lam$ with the optimal threshold level $\lam$ for each unknown vector
$\btheta$ in the class $\Theta^*_n$. This means the optimality of $\hlambda$ for the true
$\btheta$ whenever $\btheta\in \Theta^*_n$.

Property (\ref{def-adaptive-ratio-opt}) is called weak adaptive ratio optimality
when the following strong adaptive ratio optimality is also considered:
\bel{strong-ratio-opt}
\sup_{\btheta\in\Theta^*_n}\f{E_{\btheta}\|s_{\hlambda}(\bX)-\btheta\|^2}
{E_{\btheta}\inf_{\lambda\ge0}\|s_\lambda(\bX)-\btheta\|^2}\le1+o(1).
\eel
The weak and strong adaptive ratio optimality properties are of a more
appealing type than adaptive minimaxity since it applies directly to the true unknown,
instead of the worst case scenario in individual classes of unknowns.
The strong adaptive ratio optimality is even more appealing since it applies to
both the given parameter vector $\btheta$ and given data $\bX$. Define
\bel{L_0n}
L_{0,n} = (\log n)^{-3/2}\left(\delta_{1,n}(\log n)^{3/2}+(\log\log n)\sqrt{\log n}\right)^{1+\delta_{1,n}}.
\eel

\begin{thm}\label{th-adaptive-ratio-opt}
Let $\bX\sim N(\btheta,\bI_n)$ under $P_{\btheta}$ with unknown
$\btheta\in\real^n$ and $\hbtheta=s_{\hlambda}(\bX)$ be the soft threshold estimator
(\ref{hbtheta}) with a threshold level $\hlambda$ satisfying (\ref{hlam-simple}). Suppose
\bel{reg-condition-risk}
\frac{L_{0,n}}{n^{\delta_{1,n}}} \ll\inf_{\lambda\ge0}E_{\btheta}\|s_\lambda(\bX)-\btheta\|^2\ll n.
\eel
Then, $\hlambda$ approximates the optimal fixed threshold
level in the sense that
\bel{adaptive-ratio-optimality}
\limsup_{n\to\infty}\f{E_{\btheta}\|s_{\hlambda}(\bX)-\btheta\|^2}
{\inf_{\lambda\ge0}E_{\btheta}\|s_\lambda(\bX)-\btheta\|^2}\le 1.
\eel
\end{thm}

Theorem \ref{th-adaptive-ratio-opt} allows the FDR rule with $\hlam = \hxi_1$ and
$\delta_{1,n}=\delta_{2,n}=0$.  In this case, the lower risk bound for adaptive estimation is
equivalent to
\bes
\|\btheta\|^2 \gg L_{0,n} = \frac{\log\log n}{\log n}.
\ees
We prove in Lemma \ref{lm-rho_G(1)to0-part} in the next section
that (\ref{reg-condition-risk}) holds if and only if
\bel{remark-adaptive-ratio-opt-condition}
\|\btheta\|^2 \gg \frac{L_{0,n}}{n^{\delta_{1,n}}},\quad
\int_{|u|\le\eps}G_n(du)=\frac{\#\big\{i\le n\colon|\theta_i|\le\eps\big\}}{n} \to 1\
\eel
for all $\eps>0$ where $G_n$ is the nominal empirical prior as in (\ref{thresh-EB}).
In fact, we prove in Lemma \ref{lm-rho_G(1)to0-part} that for any constant $c\in (0,1)$
\bes
\sup\Big\{\Big|\frac{\inf_{\lambda\ge0}E_{\btheta}\|s_\lambda(\bX)-\btheta\|^2}{\|\btheta\|^2} - 1\Big|:
 \frac{\|\btheta\|^2}{2\log n} \le 1-c\Big\} \to 0,
\ees
so that the lower bounds in (\ref{reg-condition-risk}) and (\ref{remark-adaptive-ratio-opt-condition}) are equivalent.

\begin{example}\label{example-1} Consider the common mean model where $\theta_i=\mu$ for all $i\le n$
and $\|\btheta\|^2=n\mu^2$.
The risk of ${\bar X} = n^{-1}\sum_{i=1}^n X_i$ is
\bes
\hbox{$\hbox{\rm risk}({\bar X}) = E_{\btheta}\sum_{i=1}^n ({\bar X}-\theta_i)^2
=n E_{\btheta}({\bar X}-\mu)^2 =1.$}
\ees
When $|\mu| < (1-c)/n^{1/2}$ for a fixed $c\in (0,1)$,
the FDR soft threshold estimator $s_{\hlam}(\bX)$ outperforms ${\bar X}$ since
the risk of the FDR soft threshold estimator is no greater than
$(1+o(1))(1-c)^2 < 1= \hbox{\rm risk}({\bar X})$. Moreover, for $L_{0,n}/n^{\delta_{1,n}}\ll n\mu^2\ll 1$,
\bes
E_{\btheta}\|s_{\hlam}(\bX)-\btheta\|^2
\approx n\mu^2 \approx \inf_{\lam} n E_{\btheta}(s_\lam({\bar X})-\mu)^2.
\ees
We will follow from our general theory in the next subsection that $L_{0,n}/n^{\delta_{1,n}}\ll n\mu^2\ll 1$
also implies $E_{\btheta}\|t_{\hlam}(\bX)-\btheta\|^2 \approx n\mu^2\ll 1 = \hbox{\rm risk}({\bar X})$
for the FDR firm threshold estimator and a class of FDR smooth threshold estimators between the soft and firm.
\end{example}

Theorem \ref{th-adaptive-ratio-opt} provides the adaptive ratio
optimality in classes
\bel{class^*}
\Theta^*_n=\Big\{\btheta\in\real^n\colon
M_nL_{0,n}/n^{\delta_{1,n}} \le \inf_{\lambda\ge0}E_{\btheta}\|s_\lambda(\bX)-\btheta\|^2\le \eta_n n \Big\}
\eel
for any constant sequences satisfying $M_n\to\infty$ and $\eta_n \to 0$.
This result is a consequence of an oracle inequality in Section \ref{sec-oracle},
which uniformly bounds
\bel{regret}
\hbox{\rm regret}_{\btheta}(\hbtheta)=\f{1}{n}E_{\btheta}\|s_{\hlambda}(\bX)-\btheta\|^2
-\f{1}{n}\inf_{\lambda\ge0}E_{\btheta}\|s_\lambda(\bX)-\btheta\|^2
\eel
for the soft threshold estimator $\hbtheta=s_{\hlambda}(\bX)$ in (\ref{hbtheta}) and (\ref{hlam-simple}).

It follows from the fundamental theorem of empirical Bayes in the
compound decision theory that for any estimating function $t(x)$,
\bel{compound}
 \f{1}{n} E_{\btheta}\|t(\bX)-\btheta\|^2 = \int\int
\Big(t(x)-\theta\Big)^2\varphi(x-\theta)\,dx\,G_n(d\theta), 
\eel 
where $G_n(du) = n^{-1}\sum_{i=1}^n I\{\theta_i\in du\}$ is the unknown
nominal empirical prior \cite{Robbins51,Zhang03}. Let
$\lambda_{G_n}$ be the minimizer of (\ref{compound}) given $G_n$ for
the soft threshold estimator $t(x)=s_\lam(x)$. If $G_n$ were known,
$s_{\lambda_{G_n}}(\bX)$ could be used to achieve
$\inf_{\lambda\ge0}E_{\btheta}\|s_\lambda(\bX)-\btheta\|^2$. Thus,
as in \cite{Robbins51, Robbins56}, (\ref{regret}) can be viewed as
the regret of not knowing the nominal prior $G_n$ when one is
confined to using a soft threshold estimator.

We prove that the strong adaptive ratio optimality holds for (\ref{hbtheta}) and (\ref{hlam-simple})
in a slightly smaller range of the minimum soft threshold risk.

\begin{thm}\label{th-min-loss}
Let $\bX, \btheta, P_{\btheta}$ and $\hlambda$ be the same as in Theorem \ref{th-adaptive-ratio-opt}. Suppose
\bel{th-min-loss-1}
\log n \ll\inf_{\lambda\ge0}E_{\btheta}\|s_\lambda(\bX)-\btheta\|^2\ll n.
\eel
Then, $\hlambda$ approximates the optimal threshold level for the true $\btheta$ and almost all realizations
of data $\bX$ in the sense that
\bel{th-min-loss-2}
\lim_{n\to\infty}\f{E_{\btheta}\|s_{\hlambda}(\bX)-\btheta\|^2}
{E_{\btheta}\inf_{\lambda\ge0}\|s_\lambda(\bX)-\btheta\|^2}=1.
\eel
\end{thm}

An immediate consequence of Theorem \ref{th-min-loss} is the strong adaptive ratio optimality
in all classes in (\ref{class^*}) with $L_{0,n}/n^{\delta_{1,n}}$ replaced by $\log n$.

\subsection{Adaptive smooth threshold estimation}
Consider smooth threshold estimators $t_\lam(x)$ satisfying the following conditions:
\bel{smooth-thresh-cond}
& \begin{cases} 
s_\lam(x)\le t_\lam(x) \le h_\lam(x),& x \ge 0,
\cr h_\lam(x) \le t_\lam(x) \le s_\lam(x),& x < 0,
\end{cases}
\\ \label{smooth-cond} &\qquad \ \ \ \begin{cases} 
0\le t_\lam(y) - t_{\lam}(x)\le \kappa_0(y-x),& x<y,
\cr \big|t_\lam(x) - t_{\lam'}(x)\big|\le \kappa_1|\lam-\lam'|,&
\end{cases}
\eel
with constants $\kappa_0\in [1,2)$ and $\kappa_1 < \infty$, where
$s_\lam(x)=\sgn(x)(|x|-\lam)_+$ and $h_\lam(x)=xI\{|x|>\lam\}$ are
the soft and hard threshold estimators.
The following theorem asserts that given the unknown mean vector $\btheta$,
the risk of the FDR smooth threshold estimator,
\bel{FDR-smooth}
\hbtheta = t_{\hlam}(\bX),
\eel
with the threshold level $\hlam$ in (\ref{hlam-simple}),
is within a small fraction of the minimum risk of nonadaptive soft threshold estimator
when the risk is between $n^{-\delta_n}$ and $\delta_n n$ for any $\delta_n\to 0+$.

\begin{thm}\label{th-smooth-thresh}
Let $\bX\sim N(\btheta,\bI_n)$ under $P_{\btheta}$ with unknown
$\btheta\in\real^n$ and $\hbtheta=t_{\hlambda}(\bX)$ be the smooth threshold estimator
(\ref{FDR-smooth}) with a threshold level $\hlambda$ satisfying (\ref{hlam-simple}) and
threshold functions satisfying (\ref{smooth-thresh-cond}) and (\ref{smooth-cond}). Suppose
\bel{th-smooth-1}
\frac{L_{0,n}}{n^{\delta_{1,n}}} \ll\inf_{\lambda\ge0}E_{\btheta}\|s_\lambda(\bX)-\btheta\|^2\ll n.
\eel
Then, the FDR smooth threshold estimation is no worse than the
optimal nonadaptive soft threshold estimation in the sense that
\bel{th-smooth-2}
\limsup_{n\to\infty}\f{E_{\btheta}\|t_{\hlambda}(\bX)-\btheta\|^2}
{\inf_{\lambda\ge0}E_{\btheta}\|s_\lambda(\bX)-\btheta\|^2}\le 1.
\eel
\end{thm}

Condition (\ref{smooth-thresh-cond}) confines the estimator
$t_\lam(x)$ to the interval between the soft and hard threshold
estimators, while condition (\ref{smooth-cond}) imploses the
Lipschitz condition on $t_\lam(x)$. Both conditions hold for the
soft threshold estimator with $\kappa_0=\kappa_1=1$ and the firm
threshold estimator \cite{GaoB97} \bel{firm} f_\lam(x) =
\sgn(x)\min\Big\{|x|, \kappa_0(|x|-\lam)_+\Big\} \eel with $1 <
\kappa_0=\kappa_1< 2$. In fact conditions (\ref{smooth-thresh-cond})
and (\ref{smooth-cond}) imply that the estimator $t_\lam(x)$ must
lie between the soft and firm threshold estimators. The firm
threshold estimator can be written in the penalized form as \bes
f_\lam(x) = \inf_{\mu}\Big\{(x-\mu)^2/2 + \rho_\lam(\mu)\Big\} \ees
where $\rho_\lam(\mu) =\lam^2 \int_0^{|\mu|/\lam} (1-x/\gamma)_+dx$
is the minimax concave penalty \cite{Zhang10} with $\gamma =
1+1/\kappa_0$. The The smoothness condition (\ref{smooth-cond}) with
$c_0<2$ rules out the hard threshold estimator, which has
discontinuities at $x=\pm\lam$.

\subsection{Adaptive minimaxity with FDR smooth thresholding}\label{subsec-minimax}
For vector classes $\Theta \subset \real^n$, the minimax risk is
\bel{minimax-risk}
\scrR(\Theta)=\inf_{\bdelta}\sup_{\btheta\in\Theta}E_{\btheta}\|\bdelta(\bX)-\btheta\|^2,
\eel
where the infimum is taken over all Borel mappings $\bdelta\colon\real^n\to\real^n$.
An estimator $\hbtheta$ is asymptotically minimax with respect to a sequence of vector classes
$\Theta_n\subset\real^n$ if
\bel{def-adaptive-minimax}
\f{\sup_{\btheta\in\Theta_n}E_{\btheta}\|\hbtheta-\btheta\|^2}{\scrR(\Theta_n)} = 1+o(1).
\eel
The estimator is adaptive minimax if (\ref{def-adaptive-minimax}) holds uniformly with
a broad collection of sequences $\big\{\Theta_n\subset\real^n\big\}$ of parameter classes.

Define $\ell_p$ balls \bel{def-l_p-ball}
\Theta_{p,C,n}=\Big\{\btheta=(\theta_1,\ldots,\theta_n)\colon
\f{1}{n}\sum_{i=1}^n|\theta_i|^p\le C^p\Big\}, \eel with the
interpretation $n^{-1}\#\{i\le n:\theta_i\neq 0\}\le C$ for the
$\ell_0$ ball. The quantity $C$ is the length-normalized or
standardized radius of the $\ell_p$ ball. For $p>0$,
(\ref{def-l_p-ball}) is called the strong $\ell_p$ ball and denoted
by $\Theta_{p,C,n}^s$ when the following weak $\ell_p$ ball is also
considered as in \cite{Johnstone94,AbramovichBDJ06}:
\bel{weak-l_p-ball} \Theta_{p,C,n}^w=\Big\{\btheta\colon
|\theta_{(k)}|(k/n)^{1/p}\le C \Big\}, \eel where $|\theta_{(1)}|\ge
\ldots\ge |\theta_{(n)}|$ are the ordered absolute values of the
components of $\btheta$. We use $\Theta_{p,C,n}^{s,w}$ to denote
both strong and weak $\ell_p$ balls when a statement applies to both
types of balls.

\begin{thm}\label{th-adaptive-minimax}
Let $\bX, \btheta, P_{\btheta}$ and $\hlambda$ be the same as in Theorem \ref{th-adaptive-ratio-opt}.
Let $\scrR(\Theta)$ be the minimax risk (\ref{minimax-risk}),
$M_n'\to\infty$, $\eta_n'\to 0$, $L_{0,n}$ be as in (\ref{L_0n}), and
\bes
\Omega_n^{s,w} &=&\Big\{(p,C): 0 < p\le 2-c^{s,w},
M_n'(L_{0,n}/n^{\delta_{1,n}})^{p/2}n^{-1} \le C^{p}\le \eta_n'\Big\},
\cr \Omega_{0,n} &=& \Big\{(p,C): p=0, 1/n \le C \le \eta_n'\Big\},
\ees
with $c^{s,w}=c^{s}=0$ for strong balls and any $c^{s,w}=c^{w}\in (0,1)$ for weak balls.
Let $t_\lam(x)$ satisfy (\ref{smooth-thresh-cond}) and (\ref{smooth-cond}).
Then, for both strong and weak $\ell_p$ balls,
\bel{adaptive-minimaxity}
\lim_{n\to\infty} \sup_{(p,C)\in\Omega_n^{s,w}\cup\Omega_{0,n}}
\f{\sup_{\btheta\in\Theta_{p,C,n}^{s,w}}E_{\btheta}\|t_{\hlambda}(\bX)-\btheta\|^2}
{\scrR(\Theta_{p,C,n}^{s,w})}=1.
\eel
\end{thm}

For $0 < p\le 2-c^{s,w}$, Theorem~\ref{th-adaptive-minimax} asserts the adaptive minimaxity of the FDR
smooth threshold estimator in strong and weak $\ell_p$ balls when
\bes
(L_{0,n}/n^{1+\delta_{1,n}})^{p/2} \ll nC^{p}\ll n
\ees
with $\delta_{1,n}\to 0$ and logarithmic $L_{0,n}$ satisfying $L_{0,n}=(\log\log n)/\log n$ for $\delta_{1,n}=0$.
For $p=0$, Theorem \ref{th-adaptive-minimax} asserts the adaptive minimaxity of the FDR
smooth threshold estimator in $\ell_0$ balls when $1\le \|\btheta\|_0\ll n$.

Let $p'=pI\{p>0\}+I\{p=0\}$ and $\lam_{p,C,n}=\sqrt{2\log\big( \min(n,1/C^{p'})\big)}$.
The minimax risk for the strong and weak $\ell_p$ balls can be expressed as
\bel{minimax-risk-formula}\label{DJ-minimaxity}
\scrR(\Theta_{p,C,n}^{s,w})
&=& (1+\epsilon_{p,C,n}^{(1)})\sup\Big\{E_{\btheta}\|s_{\lam_{p,C,n}}(\bX)-\btheta\|^2: \btheta\in\Theta_{p,C,n}^{s,w}\Big\}
\cr &=& (1+\epsilon_{p,C,n}^{(2)})\sup\Big\{\|\btheta\|^2: \btheta\in\Theta_{p,C,n}^{s,w}, \|\btheta\|_\infty\le \lam_{p,C,n}\Big\}
\\ \nonumber &=& (1+\epsilon_{p,C,n}^{(3)})\inf_\lam\sup\Big\{E_{\btheta}\|h_{\lam}(\bX)-\btheta\|^2:
\btheta\in\Theta_{p,C,n}^{s,w}\Big\} \eel such that for the
$c^{s,w}$ and $\eta_n'$ in Theorem \ref{th-adaptive-minimax} \bes
\min_{n\to\infty}
\max_{k=1,2,3}\sup\Big\{\big|\epsilon_{p,C,n}^{(k)}\big|: 0\le p\le
2-c^{s,w}, C^{p'}\le \eta_n'\Big\} = 0, \ees where
$h_\lam(x)=xI\{|x|\ge \lam\}$ is the hard threshold estimator
\cite{DonohoJ94,Johnstone94,Zhang12}. Since
$\inf_\lam\sup_{\btheta\in\Theta} \ge
\sup_{\btheta\in\Theta}\inf_\lam$ for any risk function,
Theorem~\ref{th-adaptive-minimax} is almost a direct consequence of
Theorem \ref{th-adaptive-ratio-opt} and the first part of
(\ref{minimax-risk-formula}).

Statement (\ref{minimax-risk-formula}) asserts that the minimax risk
is uniformly approximately attained by either the soft  or hard
threshold estimator and that the least favorable configuration in
$\Theta_{p,C,n}^{s,w}$ is approximately attained when individual
$|\theta_i|$ concentrate at the largest possible values below a
nearly optimal threshold level. The results in
(\ref{minimax-risk-formula}) were proved in \cite{DonohoJ94} for
strong balls with $C^{p'}\gg (\log n)^{p'/2}/n$, in
\cite{Johnstone94} for weak balls with $C^{p'}\gg (\log
n)^{p'/2}/n$, and in \cite{Zhang12} for both strong and weak balls
with $C^{p'} = O(1) (\log n)^{p'/2}/n$ and explicitly stated
uniformity for the entire range of $(p,C)$.

The second part of (\ref{minimax-risk-formula}) provides an approximate formula
for the minimax risk in terms of $(p,C,n)$.
For $C^{p'}\gg (\log n)^{p/2}/n$, the formula can be more explicitly written as
\bel{minimax-0}
\scrR(\Theta_{p,C,n}^{s,w}) = (1+\epsilon_{p,C,n})M^{s,w}_p nC^{p'}\lam_{p,C,n}^{2-p},
\eel
where $M^{s,w}_p=1$ for strong balls and $M^{s,w}_p=2/(2-p)$ for weak balls.

As mentioned in the introduction, adaptive minimax estimation of
normal means in $\ell_p$ balls with $0\le p<2$ have been considered
in \cite{DonohoJ95,BirgeM01,JohnstoneS04,Zhang05,AbramovichBDJ06,
JiangZ09,WuZ13} and many others. The following results are most
closely related to Theorem~\ref{th-adaptive-minimax}: adaptive
minimaxity of the FDR hard threshold estimator for $(\log
n)^\gamma/n \le C^{p'}\le n^{-\kappa}$ with $\kappa>0$, $\gamma = 5$
in \cite{AbramovichBDJ06}, and $\gamma=4.5$ in \cite{WuZ13};
adaptive minimaxity of the GMLEB for $(\log n)^{4+p/2+3/p}/n \ll
C^{p}\ll n^{p}/(\log n)^{4+9p/2}$ for $p>0$ in \cite{JiangZ09};
adaptive rate minimaxity of the generalized $C_p$ for $1/n\le
O(1)C^{p'}$ in \cite{BirgeM01}; adaptive rate minimaxity of the
EBThresh for $(\log n)^2/n \le O(1)C^{p'}$ and a modified EBThresh
for $(\log n)^{p'/2}/n = O(1)C^{p'}$ in \cite{JohnstoneS04}. The
uniformity in $(p,C)$ of the results in
\cite{JohnstoneS04,AbramovichBDJ06,JiangZ09,WuZ13} seems to follow
from (\ref{minimax-risk-formula}) and their proofs, possibly with
some careful modification. It follows from (\ref{minimax-0}) that
the ranges of $C^{p'}$ here and those of the risk given in the
introduction are equivalent for the respective cited results.

\def\theequation{3.\arabic{equation}}
\setcounter{equation}{0}
\section{Analysis of the FDR smooth threshold estimator}\label{analysis}
As mentioned in the introduction, the compound estimation of normal means is
closely related to the Bayes estimation of a single normal mean. Let $G$ be
a prior distribution. In the Bayes problem, we
estimate a univariate random parameter $\theta$ based on a univariate
observation $X$ such that
\bes X|\theta\sim N(\theta,1),\quad \theta\sim G.
\ees

The Bayes risk of the soft threshold estimator $s_\lambda(X)$,
with fixed $\lambda$, is
\bel{Bayes-risk} 
R_G(\lambda)=\int R({u},\lambda)G(d{u}),
\eel
where $R(\mu,\lambda)$ is the conditional risk given $\theta=\mu$,
\bel{risk-soft-threshold}
&& R(\mu,\lambda)=E_{\mu}(s_\lambda(N(\mu,1))-\mu)^2=\int\Big(s_\lambda(x+\mu)-\mu\Big)^2\varphi(x)dx.
\eel
The nominal empirical Bayes prior, which naturally matches the unknown mean vector
$\btheta=(\theta_1,\ldots,\theta_n)$, is defined as
\bel{G_n}
G_n(t)=\f{1}{n}\sum_{i=1}^n I\big\{\theta_i\le t\big\}.
\eel
With the above notation, (\ref{compound}) with $t(x)=s_\lam(x)$ can be written as
\bel{EB}
R_{G_n}(\lambda)=\f{1}{n}E_{\btheta}\|s_\lambda(\bX)-\btheta\|^2.
\eel
We denote the Bayes optimal soft threshold risk and level for prior $G$ by
\bel{lam_opt}
\eta_G = \min_{\lam} R_{G}(\lam),\quad \lambda_{G} = \argmin_\lam R_{G}(\lam).
\eel
It follows immediately from (\ref{EB}) that $\lam_{G_n}$ is the optimal deterministic
soft threshold level when $\btheta$ is the true mean vector and
\bel{target}
\eta_{G_n} = R_{G_n}(\lambda_{G_n}) = \f{1}{n}\,\inf_{\lam\ge 0}E_{\btheta}\|s_\lambda(\bX)-\btheta\|^2.
\eel
For smooth threshold functions satisfying (\ref{smooth-thresh-cond}) and (\ref{smooth-cond}), we define
\bel{R-sm}
R^{(sm)}_{G}(\lam)=\int \int \Big(t_\lam(x+u)-u\Big)^2\varphi(x)\, dx \,G(du).
\eel

Our analysis requires a concentration inequality to bound from the above the difference
$\|t_{\hlam}(\bX)-\btheta\|^2/n - R^{(sm)}_{G_n}(\hlam)$
and exponential inequalities to bound from the above
$R^{(sm)}_{G_n}(\hlam) - R_{G_n}(\lambda_{G_n})$. Define
\bel{R_1n}
R_{1,n}(\btheta,\lam) = E_{\btheta}\| t_{\hlam}(\bX)-\btheta\|/\sqrt{n}.
\eel
The concentration inequality actually bounds the
difference between random variables $\|t_{\hlam}(\bX)-\btheta\|/\sqrt{n}$ and
$R_{1,n}(\btheta,\hlam)$ based on the fact that for each deterministic $\lam$,
the Lipschitz norm of $\|t_{\lam}(\bX)-\btheta\|$ is no greater than $\kappa_0$ as a function of the error $\bX-\btheta$.
Since $\hlam$ is bounded by functions of the FDR rules $\{\hxi_1,\hxi_2\}$ in (\ref{xi}) and
the FDR rules are defined through the counting process $N(t)$ in (\ref{N(t)}),
exponential inequalities for $N(t)$ are used to bound the
difference between $R^{(sm)}_{G_n}(\hlam)$ and $R_{G_n}(\lambda_{G_n})$.

We provide below preliminary analysis of the Bayes risk function $R_G(\lam)$ for the soft threshold estimator,
that of $R^{(sm)}_{G_n}(\lam)$ for the smooth threshold estimator,
that of the FDR rules $\{\hxi_1,\hxi_2\}$, and the concentration inequality.
Throughout the analysis, we denote by $M^*$ a positive numerical
constant which may take different values in different appearances.

\subsection{Risk properties of soft thresholding at fixed level}
\label{risk-property-soft-threshold}
With a numerical constant ${B}_0\ge 4$, let
\bel{r_G}
&\rho_G(\lambda)=\int(u^2\wedge\lambda^2)G(du),\quad
r_G(\lambda)\equiv\rho_G(\lambda)+{B}_0\Phi(-\lambda).
\eel
We carry out an analysis of the Bayes risk $R_G(\lam)$ by studying the relationship between $R_G(\lam)$
and the more explicit $r_G(\lam)$. Parallel to (\ref{lam_opt}), we define
\bel{lam-star}
\eta_G^* = \min_{\lam} r_{G}(\lam),\quad \lambda_{G}^* = \argmin_\lam r_{G}(\lam).
\eel

\begin{lem}\label{lm-Bayes-approx}
Let $R_G(\lambda)$, $R(\mu,\lam)$, $\eta_G$, $\lambda_{G}$, $r_G(\lambda)$, $\rho_G(\lam)$,
$\eta_G^*$ and $\lam_G^*$ be as in (\ref{Bayes-risk}), (\ref{risk-soft-threshold}), (\ref{lam_opt}), (\ref{r_G})
and (\ref{lam-star}) respectively.

(i) The soft threshold risk $R(0,\lam)$ at $\mu=0$ is decreasing in $\lam$ with
\bel{R(0,lambda)-bounds}
\f{4\Phi(-\lambda)}{\lambda^2+5}\le
R(0,\lambda)\le\f{4\Phi(-\lambda)}{\lambda^2+2}.
\eel
The Bayes risk $R_G(\lambda)$ is bounded by
\bel{R_G-upper-bound}\label{r_G-bd}
R_G(\lambda)\le \rho_G\big(\sqrt{\lambda^2+1}\big)+R(0,\lambda)
\le \rho_G(\lam)+r_G(\lam)/(\lam^2\vee 1).
\eel
Consequently, the minimum Bayes risk is bounded by
\bel{eta_G-bd}
\eta_G \le \Big(1+\frac{1}{(\lam^*_G)^2\vee 1}\Big)\eta_G^*.
\eel

(ii) There exists a constant $M^*_0$ depending on $B_0$ only such that
\bel{Bayes-approx-Lambda}
\max\{\rho_G(1),\eta^*_G\}\le \Big(1 + \frac{\sqrt{8\log(\lam_G\vee e)}+M^*_0}{\lam_G\vee e}\Big)\eta_G.
\eel
\end{lem}

Lemma \ref{lm-Bayes-approx} provides an approximation of the optimal risk
$\eta_G\approx \eta^*_G$ for large $\lam_G\wedge\lam^*_G$.
We now consider small Bayes soft threshold risk.

\begin{lem}\label{lm-rho_G(1)to0-part}
Let $R_G(\lambda)$, $G_n$, $\eta_G$, $\lambda_{G}$, $r_G(\lambda)$, $\rho_G(\lam)$, $\eta_G^*$ and $\lam_G^*$
be as in (\ref{Bayes-risk}), (\ref{G_n}) , (\ref{lam_opt}), (\ref{r_G}) and (\ref{lam-star}) respectively.

(i) The quantity $\lam^*_G$ is an increasing function of $B_0$.
There exists a numerical constant $M^*$ such that
\bel{lm-rho_G(1)to0-part-1}
 \sqrt{2\log(1/\eta_n)} &\le& \left(1+ \frac{\log\log(1/\eta_n)-\log(8/7)}{4\log(1/\eta_n)}\right)\lam^*_G,
\cr \sqrt{2\log(1/\eta_n)} &\le& \left(1+ \frac{3\log\log(1/\eta_n)+\log 4}{4\log(1/\eta_n)}\right)\lam_G,
\eel
whenever $\min(\eta_G,\eta^*_G)\le\eta_n \le 1/M^*$. Moreover, under the same condition
\bel{lm-R_G-r_G-bd}
\max\Big\{\eta_G,\rho_G(\lam^*_G)\Big\} &\le& \left(1+ \frac{1}{\log(1/\eta_n)}\right) \eta_G^*
\cr &\le& \left(1+ \frac{2\sqrt{\log\log(1/\eta_n)}+M^*_0}{\sqrt{\log(1/\eta_n)}}\right)\eta_G
\eel
with a constant $M^*_0$ depending on $B_0$ only.

(ii) If $\sqrt{2\log(1/\eta_n)}\ge {B}_0/\sqrt{2\pi}$ and $\rho_G(1)\le\eta_n$, then
\bel{lm-rho_G(1)to0-part-2}
\max\Big(\eta_G-\eta_n, \eta_G^*\Big) \le
\eta_n\Big(1+2\log(1/\eta_n)\Big).
\eel

(iii) Let $z_n>0$ satisfy $z_n^{-2}\Phi(-z_n) = 1/(4n)$. Then, $z_n^2>\log n$ for $n\ge 2$,
$z_n^2 > 2\log\big(4n/\sqrt{2\pi}\big) - 3\log\big(2\log\big(4n/\sqrt{2\pi}\big)\big)$ for $n\ge 7$,
and
\bel{lm-2-3}
\lam^*_{G_n}=\infty\ \hbox{ and }\
\eta_{G_n}^* = \|\btheta\|^2/n
\eel
for $\eta^*_{G_n}\le z_n^2/n$.
Moreover, there exists a numerical positive integer $n^*$ such that (\ref{lm-2-3}) holds
whenever $\sqrt{n\eta_{G_n}} \le \sqrt{2\log n}-2\sqrt{2\log\log n}$ and $n\ge n^*$.
\end{lem}

Lemma \ref{lm-rho_G(1)to0-part} (i) provides the approximation of the optimal nonadaptive soft threshold 
risk, $\eta_G\approx \eta^*_G$,
when $\min(\eta_G,\eta^*_G)$ is small, compared with the less explicit condition of having large
$\lam_G\wedge\lam^*_G$ in Lemma \ref{lm-Bayes-approx} (ii).
Lemma \ref{lm-rho_G(1)to0-part} (ii) and Lemma \ref{lm-Bayes-approx} (ii) imply the
equivalence of the condition $\min(\eta_G,\eta^*_G)\to 0$ and the even more explicit $\rho_G(1)\to 0$,
which implies the equivalence of the upper risk bound conditions in (\ref{reg-condition-risk})
and (\ref{remark-adaptive-ratio-opt-condition}) in Proposition \ref{prop-1} (i) below.
Lemma \ref{lm-rho_G(1)to0-part} (iii) gives explicit expression of $\eta^*_{G_n}$ when
the optimal risk $\eta_{G_n}$ is smaller than a critical risk level near $2(\log n)/n$,
which implies the equivalence of the lower risk bound conditions in (\ref{reg-condition-risk})
and (\ref{remark-adaptive-ratio-opt-condition}) as described in Proposition \ref{prop-1} (ii).
Define $\Gbar(t) = \int_{|u|>t}G(du)$. For $0<\epsilon\le 1$, we have
\bes
\epsilon^2\Gbar(\epsilon) \le \rho_G(1) \le \epsilon^2+\Gbar(\epsilon).
\ees

\begin{proposition}\label{prop-1}
Let $R_G(\lambda)$, $\eta_G=\inf_\lambda R_{G}(\lambda)$,
$r_G(\lambda)$, $\rho_G(\lam)$ and $\eta_G^*=\inf_\lambda r_{G}(\lambda)$
be as in (\ref{Bayes-risk}), (\ref{lam_opt}), (\ref{r_G}) and (\ref{lam-star}) respectively.
Let $G_n$ be the nominal empirical prior in (\ref{G_n}) for the unknown vector $\btheta$.

(i) For fixed $B_0$ the following conditions are equivalent to each other:
(a) $\eta_{G}\to 0$;
(b) $\eta^*_{G}\to0$;
(c) $\rho_{G}(1)\to0$; and
(d) $\Gbar(\epsilon)\to 0$ for all $\eps> 0$.

(ii) There exists a numerical positive integer $n^*$ such that
\bes
\eta_{G_n}^* =\|\btheta\|^2/n = (1+o(1))\eta_{G_n}
\ees
whenever $\min(\eta_{G_n},\eta_{G_n}^*)\le (\sqrt{2\log n}-3\sqrt{2\log\log n})_+^2/n$ and $n\ge n^*$.
\end{proposition}

We omit the proof of Proposition \ref{prop-1} since it is a direct consequence of
Lemmas \ref{lm-Bayes-approx} and \ref{lm-rho_G(1)to0-part} as discussed above its statement.

\subsection{Risk properties of smooth thresholding at fixed level}
Let $t_\lam(x)$ be threshold functions between the soft and firm threshold estimators:
\bel{thresh-cond}
\begin{cases}
s_\lam(x) \le t_\lam(x) \le f_\lam(x),& x \ge 0
\cr f_\lam(x) \le t_\lam(x) \le s_\lam(x),& x < 0
\end{cases}
\eel
where $s_\lam(x)=\sgn(x)(|x|-\lam)_+$ and $f_\lam(x)$ is as in (\ref{firm}) with $\kappa_0\in [1,2)$.

\begin{lem}\label{lm-firm}
Suppose (\ref{thresh-cond}) holds with $\kappa_0\in [1,2)$. Let $C_0=\kappa_0/(2-\kappa_0)$. \\
(i) For all $\mu$ and $\lam\ge 0$,
\bes
|t_\lam(x)-\mu| \le \max\Big(|\mu|, C_0(|x-\mu|-\lam)_+\Big).
\ees
(ii) Suppose $EX=\mu$. Then, for all $\lam\ge 0$,
\bes
E\Big(t_\lam(X)-\mu\Big)^2\le \lam^2+2\,\Var(X).
\ees
(iii) Let $1\le q\le 2$. Suppose $EX=\mu$. Then, for all $\lam\ge 0$,
\bes
&& E\Big|t_\lam(X)-\mu\Big|^q\le \min\Big(|\mu|^q+C_0^qE(|X-\mu|-\lam)_+^q,\big(\lam^2+2\Var(X))^{q/2}\Big).
\ees
(iv) Let $\bX\sim N(\btheta,\bI)$ under $P_{\btheta}$ and $R(\mu,\lam)$, $G_n$,
$r_G(\lambda)$, $\rho_G(\lam)$ be as in (\ref{risk-soft-threshold}), (\ref{G_n}) and (\ref{r_G}) respectively
with $B_0\ge 4\vee (2C_0^2)$. Then,
\bes
n^{-1} E_{\btheta}\|t_\lam(\bX)-\btheta\|_2^2
&\le& \rho_{G_n}(\sqrt{\lam^2+2})+C_0^2 R(0,\lam)
\cr &\le& \rho_{G_n}(\lam) + 2 r_{G_n}(\lam)/(\lam^2\vee 1).
\ees
\end{lem}

It follows from Lemma \ref{lm-firm} (iv) that for the optimal $\lam^*_{G}$ and $\eta^*_{G}$ in (\ref{lam-star}),
\bes
n^{-1} E_{\btheta}\|t_{\lam^*_{G_n}}(\bX)-\btheta\|_2^2 \le \Big(1+2/(\lam^*_{G_n}\vee 1)^2\Big)\eta^*_{G_n}.
\ees
Thus, as in Lemma \ref{lm-rho_G(1)to0-part}, condition $\min(\eta_G,\eta^*_G)\le\eta_n \le 1/M^*$ implies
\bes
&& E_{\btheta}\|t_{\lam^*_{G_n}}(\bX)-\btheta\|_2^2
\cr &\le& \left(1+ \frac{2\sqrt{\log\log(1/\eta_n)}+M^*_0}{\sqrt{\log(1/\eta_n)}}\right)
\inf_\lam E_{\btheta}\|s_\lam(\bX)-\btheta\|_2^2.
\ees
This asserts that at a proper threshold level,
the risk of smooth thresholding satisfying (\ref{thresh-cond}) can not be significantly larger
than that of the optimal soft thresholding. The reverse is not true in view of the following example.

\begin{example} Let $\#\{i:\theta_i=0\}=n-1$ and $\#\{i: \theta_i = \mu\}=1$ with $\mu=4\sqrt{2\log n}$.
Lemma \ref{lm-rho_G(1)to0-part} yields
$n\eta_{G_n}^*\ge n\rho_n(\lam_{G_n}^*)\ge (\lam_{G_n}^*)^2\wedge \mu^2$, so that
\bes
n\eta_{G_n}= \inf_\lam E_{\btheta}\|s_\lam(\bX)-\btheta\|_2^2 \ge (1+o(1))2\log n.
\ees
On the other hand, for the firm threshold estimation (\ref{firm}) with $\kappa_0=3/2$ and $\lam=\sqrt{2\log n}$,
we have $C_0=\kappa_0/(2-\kappa_0)=3$ and
\bes
E_{\btheta}\|f_\lam(\bX)-\btheta\|_2^2\le 9n R(0,\lam)+ E\Big(f_\lam(N(\mu,1))-\mu\Big)^2 = O(1).
\ees
\end{example}

\subsection{Analysis of the FDR threshold level}
We discuss the relationship between the FDR threshold levels (\ref{xi})
and their population version.

A population version of the FDR can be defined as
\bes
\hbox{FDR}_{\mathrm{pop}} = \frac{E\#\{\hbox{\,falsely rejected hypotheses\,}\}}{E \#\{\hbox{\,rejected hypothesis\,}\}}.
\ees
Let $G_n$ be the nominal empirical prior in (\ref{G_n}). Define
\bel{S_G-Gbar}
S_G(t)=\int
P\big\{|N({u},1)|>t\big\}G(d{u})
\eel
for any probability distribution $G$.
If $\btheta$ has $n_0$ zero components and $H_i:\theta_i=0$ is tested by thresholding $|X_i|$
at level $t$, the population FDR is
\bes
\hbox{FDR}_{\mathrm{pop}}(t) = \frac{n_0 P\{|N(0,1)|>t\}}{\sum_{i=1}^n P\{|N(\theta_i,1)|>t\}}
= \frac{n_0 2\Phi(-t)}{nS_{G_n}(t)}.
\ees
We call $2\Phi(-t)/S_{G_n}(t)$ the nominal FDR function as its sample version.

Since this paper is concerned with estimation, the $\ell_0$ sparsity of $\btheta$ is
covered but not assumed.
Actually, we allow $n_0=\#\{i\le n: \theta_i=0\}=0$. Still, the nominal FDR function
$2\Phi(-t)/S_{G_n}(t)$ plays a crucial role in studying the FDR threshold level (\ref{xi}).
Given two nominal FDR levels $\alpha_1'$ and $\alpha_2'$, the population version of
the threshold levels (\ref{xi}) is
\bel{xi_*}
&& \xi_{1,*}=\inf\Big\{t: \frac{2\Phi(-t)}{S_{G_n}(t)} \le \alpha_1'\Big\}, \
\xi_{2,*}=\sup\Big\{t: \frac{2\Phi(-t)}{S_{G_n}(t)} \ge \alpha_2'\Big\}
\eel
with the $G_n$ in (\ref{G_n}).
We consider fixed $0<\alpha'_2<\alpha_2\le\alpha_1<\alpha'_1<1$. Let
\bel{Gbar}
\Gbar(t)=\int_{|{u}|>t}G(d{u}).
\eel
As we have mentioned in the introduction, we will present an oracle inequality
in Section \ref{sec-oracle} for a more general class of threshold rules.
This class involves certain functions $g_{1,n}$ satisfying
\bel{g_1-cond}
&& 0\le g_{1,n}(x) \le x,\ R(0,g_{1,n}(x)) \le 4\Phi(-x),\quad \forall\ x >0.
\eel

\begin{lem}\label{lm-Gauss-mix}
Let $S_G(t)$, $\Gbar(t)$, $\rho_G(t)$, $G_n$, $\xi_{1,*}$ $\xi_{2,*}$ be as in
(\ref{S_G-Gbar}), (\ref{Gbar}), (\ref{r_G}), (\ref{G_n}) and (\ref{xi_*}) respectively.

(i) Suppose $S_G(t)\neq 2\Phi(-t)$ for some $t>0$, i.e. $\Gbar(t)>0$ for some $t>0$.
Then, the nominal population FDR level, $2\Phi(-t)/S_G(t)$, is strictly decreasing in $t$
from 1 at $t=0$ to 0 as $t\to\infty$, and that for all $t>0$,
\bel{Gauss-mix-1}
\Gbar(t)/2 \le S_G(t)\le 2\Phi(-t) + \rho_G(1).
\eel
Moreover, $\lam^*_{G_n} > \xi_{2,*}$ when ${B}_0\ge 8/\alpha_2'$.

(ii) For $j=1,2$, $2\Phi(-\xi_{j,*}) \le \rho_{G_n}(1)\alpha_j'/(1-\alpha_j')$.
If (\ref{g_1-cond}) holds, then
\bel{Gauss-mix-2}
\alpha_j'S_{G_n}(g_{1,n}(t)) \le (5+t^2)2\Phi(-t),\quad \forall t\le \xi_{j,*}
\eel

(iii) Let $A_1\ge 1$ and $\theta_{*,n} =\max\{t: A_1^{1/2}\xi_{1,1}t+t^2/2\le \beta_0 \log n\}$
with a certain $\beta_0 \in (0,1/2]$ and the $\xi_{1,1}$ in (\ref{xi_jk}).  Suppose
\bel{cond-xi_*}
\eta_{G_n}^* \le \theta_{*,n}^2/n,\ n\ge 2,\ 1+n^{\beta_0-1/2}/2 \le 1/\alpha_1'.
\eel
Then, $\xi_{1,*} \ge A_1^{1/2}\xi_{1,1}$ with the population FDR threshold level (\ref{xi_*}).
Moreover, for all $t$ satisfying $2\theta_{*,n} \le t\le A_1^{1/2}\xi_{1,1}$,
\bel{Gauss-mix-3}
&& \int_t^\infty S_{G_n}^{1/2}(x)dx \le t^{-1}\sqrt{2\Phi(-t)}\Big(2+4\sqrt{n^{-1/2}+n^{\beta_0-1/2}}\Big).
\eel
\end{lem}

Lemma \ref{lm-Gauss-mix} (i) provides the monotonicity of the population FDR
as a function of the threshold level $t$ and lower and upper bounds for the population
rejection probability $S_G(t)$.
Lemma \ref{lm-Gauss-mix} (ii) provides a lower bound for the population FDR threshold level
$\xi_{j,*}$ and an upper bound for the population rejection probability at level $g_{1,n}(t)$.
Lemma \ref{lm-Gauss-mix} (iii) provides a condition under which the population FDR threshold
level is greater than the highest possible sample FDR threshold level $\xi_{1,1} = - \Phi^{-1}(\alpha_1/n)$
at the nominal FDR level $\alpha_1$.
We note that the third condition in
(\ref{cond-xi_*}) holds for all $n$ and $\beta_0\le 1/2$ when $\alpha_1'\le 2/3$.

Since $0<\alpha'_2<\alpha_2\le\alpha_1<\alpha'_1<1$ and smaller false discovery error
requires higher threshold level, we expect
$\hxi_1 \le \xi_{1,*}\le \xi_{2,*} \le\hxi_2$ with large probability.
This is verified in the following lemma.

\begin{lem}\label{lm-exp-inq}
Let $\bX\sim N(\btheta,\bI_n)$ under $P_{\btheta}$. Let $N(t)$ be as
defined in (\ref{N(t)}). Let $\xi_{j,k}$, $\hxi_j$ and $\xi_{j,*}$ be as
defined in (\ref{xi_jk}), (\ref{xi}) and (\ref{xi_*}). Then

(i) For all $\xi_k>0$ (e.g. $\xi_k=\xi_{1,k}$ or \,$\xi_k=\xi_{2,k}$),
\bel{exp-inq-1}
&P_{\btheta}\Big\{\sgn(k-E_{\btheta}N(\xi_k))(N(\xi_k)-k)\ge0\Big\}
\le\exp(-\nu_k k),\eel where
$\nu_k=E_{\btheta}N(\xi_k)/k-1-\log(E_{\btheta}N(\xi_k)/k)>0$.

(ii) Let $\nu_{1,*}=\alpha_1/\alpha'_1-1-\log(\alpha_1/\alpha'_1)$. For all $\xi_{1,k}\le\xi_{1,*}$,
\bel{exp-inq-2}
P_{\btheta}\Big\{ \hxi_1\le \xi_{1,k}\Big\} \le P_{\btheta}\Big\{ N(\xi_{1,k})\ge k\Big\} \le \exp(-\nu_{1,*}k).
\eel
Let $\nu_{2,*}=\alpha_2/\alpha'_2-1-\log(\alpha_2/\alpha'_2)$. For all $\xi_{2,k}\ge\xi_{2,*}$,
\bel{exp-inq-3}
P_{\btheta}\Big\{ \hxi_2 \ge \xi_{2,k} \Big\}
\le P_{\btheta}\Big\{ N(\xi_{2,k})\le k \Big\}\le \exp(-\nu_{2,*}k).
\eel
\end{lem}

\subsection{Gaussian isoperimetric inequality}
Here we provide large deviation bounds, based on the Gaussian
isoperimetric inequality \cite{Borell75, vanderVarrtW96,
MassartP07}, for the difference between the loss and risk functions
of the smooth threshold estimators satisfying
(\ref{smooth-thresh-cond}) and (\ref{smooth-cond}) at an arbitrary
random threshold level.

For real-valued functions $f$ on $\real^n$, the Lipschitz norm is defined as
\bes
\|f\|_{\mathrm{Lip}}=\sup_{\bu\neq\bv}\f{|f(\bu)-f(\bv)|}{\|\bu-\bv\|}.
\ees
The Gaussian isoperimetric inequality asserts that for
$\bZ \sim N(0,\bI_n)$ and functions $f\colon\real^n\to\real$ with $\|f\|_{\mathrm{Lip}}\le 1$,
\bel{Gaussian-iso-inq}
P\big\{|f(\bZ)-Ef(\bZ)|>t\big\} \le 2e^{-t^2/2}.
\eel
By (\ref{smooth-cond}), the smooth threshold estimator $t_\lam(x)$ has Lipschitz norm $\kappa_0$
as a function of $x$, so that the $\ell_2$ norm of the loss $\|t_\lam(\bX)-\btheta\|$ also has
Lipschitz norm $\kappa_0$ as a real valued function of $\bZ=\bX-\btheta$. This leads to the following lemma.
Let $L(\bX,\btheta,\lambda)=\|t_\lambda(\bX)-\btheta\|/\sqrt{n}$ and define
\bel{H}
\qquad R_{1,n}(\btheta,\lambda) &=&E_{\btheta}L(\bX,\btheta,\lambda),
\cr H_{\btheta}(a,b) &=&\Big|L(\bX,\btheta,Y)-R_{1,n}(\btheta,Y)\Big|^2I\big\{a\le Y<b\big\},
\eel
with a nonnegative random variable $Y$, and
\bel{H-star}
H^*_{\btheta}(\lam) = \max_{y\ge\lam}\Big|L(\bX,\btheta,y)-\|\btheta\|/\sqrt{n}\Big|^2.
\eel

\begin{lem}\label{lm-Gaussian-iso-inq2}
Let $\bX\sim N(\btheta,\bI_n)$ under $P_{\btheta}$ and
$\{S_{G}(t), G_n\}$ be as in (\ref{S_G-Gbar}) and (\ref{G_n}).
Then, for all $a=c_0<\cdots<c_m=b$ and $\pi_j \ge
P_{\btheta}\big\{c_{j-1}\le Y<c_j\big\}$,
\bel{lm-Gaussian-iso-discrete}
\sqrt{E_{\btheta} H_{\btheta}(a,b)} &\le &
2\kappa_0\bigg\{\f{2}{n}\sum_{j=1}^m \pi_j \log\Big(\frac{2e}{\pi_j\wedge 1}\Big)\bigg\}^{1/2}
\cr  && +2\kappa_1\bigg\{\sum_{j=1}^m\pi_j(c_j-c_{j-1})^2S_{G_n}(c_{j-1})\bigg\}^{1/2}.
\eel
In particular, for all integers $m\ge1$,
\bel{lm-Gaussian-iso-discrete-coll}
\sqrt{E_{\btheta}H_{\btheta}(a,b)}
\le 2\kappa_0\bigg\{\f{2}{n}\log(2em)\bigg\}^{1/2}
+2\kappa_1\bigg\{\f{(b-a)^2}{m^2}S_{G_n}(a)\bigg\}^{1/2}.
\eel
Moreover, for all $\lam>0$,
\bel{lm-Gaussian-iso-star}
 \sqrt{E_{\btheta}H^*_{\btheta}(\lam)} \le \kappa_1\int_\lam^\infty S_{G_n}^{1/2}(t)dt
\eel
and with the $R^{(sm)}_{G}(\lam)$ in (\ref{R-sm})
\bel{lm-Gaussian-iso-basic}
R_{1,n}^2(\btheta,\lambda)\le E_{\btheta}L^2(\bX,\btheta,\lambda) = R_{G_n}^{(sm)}(\lam)
\le R_{1,n}^2(\btheta,\lambda)+\f{4\kappa_0^2}{n}.
\eel
\end{lem}

\def\theequation{4.\arabic{equation}}
\setcounter{equation}{0}
\section{An oracle inequality}\label{sec-oracle}
We provide an oracle inequality for a more general class of threshold levels.
Let $g_{1,n}(x)$ be a sequence of functions satisfying the following conditions:
\bel{cond-g_1}
&& 0\le g_{1,n}(x) \le x,\  0\le (d/dx)g_{1,n}(x)\le M_0,\ \forall x >0,
\cr && R(0,g_{1,n}(x)) \le \min\left\{4\Phi(-x), \frac{M_0\Phi(-x)}{(x^{c_{1,n}}+2)(1\vee \log x)^{c_{2,n}}}\right\}
\\ \nonumber && 0 < c_{1,n}\le 2,\ |c_{2,n}|\le M_0,\ c_{2,n}\le 0\ \hbox{ for } c_{1,n}=2,
\eel
where $M_0$ is a numerical constant.
Recall that $R(0,x)=E(|N(0,1)|-x)_+^2$.
Since $R(0,x)\le 4\Phi(-x)/(x^2+2)$ for all $x>0$
by (\ref{R(0,lambda)-bounds}) of Lemma \ref{lm-Bayes-approx} (i),
(\ref{cond-g_1}) holds for $g_{1,n}(x)=x$ with $M_0=4$, $c_{1,n}=2$ and $c_{2,n}=0$.

Threshold levels of the following form will be considered:
\bel{hlam}
\sqrt{1+\delta_{1,n}}g_{1,n}(\hxi_1) \le\hlam \le \sqrt{1+\delta_{2,n}}\hxi_2
\eel
with $0\le \delta_{1,n}\le\delta_{2,n}$. Compared with (\ref{hlam-simple}),
the lower bound in (\ref{hlam}) is smaller with $g_{1,n}(\hxi_1)\le \hxi_1$ and the upper bound
in (\ref{hlam}) is larger since $\delta_{2,n}\to 0$ is no longer assumed. We note that
$\delta_{2,n}\to 0$ is necessary for attaining the optimal constant factor in our analysis but not
for rate optimality.

We prove in the Appendix that for large $x$, condition (\ref{cond-g_1}) implies
\bel{cond-g_1-a}
\Phi\Big(-\sqrt{1+\delta_{1,n}}g_{1,n}(x)\Big)
\le  x^{-1}\left(\frac{M^*x^3\Phi(-x)}{x^{c_{1,n}}(\log_+x)^{c_{2,n}}}\right)^{1+\delta_{1,n}}
\eel
with a numerical constant $M^*$ depending on $M_0$ only. Define
\bel{L_2n}
 L_{2,n} = (\log n)^{-3/2}\left(\frac{\delta_{1,n}(\log n)^{(5-c_{1,n})/2}}{(\log_+\log n)^{c_{2,n}}}
+\frac{(\log n)^{(3-c_{1,n})/2}}{(\log_+\log n)^{c_{2,n}-1}} \right)^{1+\delta_{1,n}}.
\eel

\begin{thm}\label{thm-4}
Let $\alpha'_2<\alpha_2\le\alpha_1<\alpha'_1<1$, $\beta_0\le 1/2$ and $A$ be fixed positive constants.
Let $\bX\sim N(\btheta,\bI_n)$ under $P_{\btheta}$ and $\hlam$ be a threshold level satisfying
(\ref{hlam}) with $\hxi_1$ and $\hxi_2$ in (\ref{xi}), constants $0\le \delta_{1,n}\le\delta_{2,n}\wedge A$
and functions $g_{1,n}$ satisfying (\ref{cond-g_1}).
Let $t_\lam(x)$ be threshold functions satisfying (\ref{smooth-thresh-cond}) and (\ref{smooth-cond})
with $C_0=\kappa_0/(2-\kappa_0)$,
$G_n$ be as (\ref{G_n}),
$r_G(\lam)$ as in (\ref{r_G}) with ${B}_0= (8/\alpha'_2)\vee(2C_0^2)$,
and $\eta^*_G$ as in (\ref{lam-star}).
Assume $\max_{1\le n < n_*}(1+\delta_{2,n})\le A$ with
$n_* = \min\big\{n: 1+n^{\beta_0-1/2}/2 \le 1/\alpha_1', n\ge 2\big\}$. Then,
\bel{thm-4-1}
&& \sqrt{E_{\btheta}\|t_\hlam(\bX)-\btheta\|^2/n}
\cr & \le & \sqrt{(1+ \delta_{2,n})\eta_{G_n}^*}
+ M^*\sqrt{(1+ \delta_{2,n})\tau^*_{1,n}\eta_{G_n}^*+\tau^*_{2,n}},
\eel
where $M^*$ is a constant depending on $\{\alpha_1',\alpha_2',\beta_0,A,M_0\}$ only
with the $M_0$ in (\ref{cond-g_1}),
$\tau^*_{2,n} = L_{2,n}/n^{1+\delta_{1,n}}$ and with $L_{1,n}^*=e\vee \log(1/\eta_{G_n}^*)$
\bel{tau_n}
\tau_{1,n}^* = \max\left(\frac{\log(e\vee \log n)}{1\vee \log n},
\frac{(\log L_{1,n}^*)^{-c_{2,n}}}{(L_{1,n}^*)^{c_{1,n}/2}},\frac{1}{L_{1,n}^*}\right).
\eel
\end{thm}

We note that $n_*=2$ when $\alpha_1'\le 2/3$.

Let $\eta_{G_n}=\inf_{\lam\ge 0}E_{\btheta}\|s_\lam(\bX)-\btheta\|^2/n$
defined through (\ref{EB}), (\ref{lam_opt}) and (\ref{G_n})
and $\eta_{G_n}^*=\inf_\lam r_{G_n}(\lam)$ defined through (\ref{r_G}) and (\ref{lam-star}).
It follows from (\ref{lm-R_G-r_G-bd}) of Lemma \ref{lm-rho_G(1)to0-part} that
$\eta_{G_n}$ and $\eta_{G_n}^*$ are within
a small fraction of each other when $\eta_{G_n}\wedge\eta_{G_n}^*$ is small,
so that Theorem \ref{thm-4} with $t_\lam(x)=s_\lam(x)$ implies adaptive ratio optimality
and minimaxity of the FDR soft threshold estimator when $\eta_{G_n}\vee\delta_{2,n} \to 0$ 
and $\tau^*_{2,n}\ll \eta_n$.
The following corollaries provide a more general and more explicit version of
Theorems \ref{th-adaptive-ratio-opt}, \ref{th-smooth-thresh} and \ref{th-adaptive-minimax}.

\begin{cor}\label{cor-1} Let $\bX\sim N(\btheta,\bI_n)$ under $P_{\btheta}$ and
$\hlam$ be a threshold level satisfying (\ref{hlam}) with $0\le\delta_{1,n}\le\delta_{2,n}\to 0$
and functions $g_{1,n}$ satisfying (\ref{cond-g_1}). Let
$\eta_n\ge \eta_{G_n}= R_{G_n}(\lam_{G_n})=\min_{\lam\ge 0}E_{\btheta}\|s_\lam(\bX)-\btheta\|^2/n$. Let
\bes
\tau_{1,n} = \max\left(\delta_{2,n},\frac{\log(e\vee \log n)}{1\vee \log n},
\frac{(\log L_{1,n})^{-c_{2,n}}}{L_{1,n}^{c_{1,n}/2}},\frac{\sqrt{\log L_{1,n}}}{\sqrt{L_{1,n}}}\right)
\ees
with $L_{1,n}=e\vee \log(1/\eta_n)$ and
$\tau^*_{2,n} = L_{2,n}/n^{1+\delta_{1,n}}$ with $L_{2,n}$ in (\ref{L_2n}).
Let $t_\lam(x)$ be functions satisfying (\ref{smooth-thresh-cond}) and (\ref{smooth-cond})
with $C_0=\kappa_0/(2-\kappa_0)$.
Then,
\bel{cor-1-1}
\sqrt{E_{\btheta}\|t_\hlam(\bX)-\btheta\|^2/n}
\le \sqrt{\eta_{G_n}} + M^*\sqrt{\tau_{1,n}\eta_{G_n} + \tau^*_{2,n}},
\eel
where $M^*$ is a constant depending on $\{\alpha_1',\alpha_2',\beta_0,A,M_0\}$ only.
Consequently, the adaptive ratio optimality (\ref{def-adaptive-ratio-opt}) for
\bes
\Theta_n^*=\Big\{\btheta: M_nL_{2,n}/n^{1+\delta_{1,n}} \le R_{G_n}(\lam_{G_n}) \le\eta_n\Big\}
\ees
as long as $M_n\to\infty$ and $\eta_n\to 0$.
\end{cor}

\begin{cor}\label{cor-2} Let $\bX$, $\btheta$, $P_{\btheta}$, $t_\lam(x)$ and
$\hlam$ be as in Corollary \ref{cor-1}.
Then, the conclusion of Theorem \ref{th-adaptive-minimax} holds
for $\hbtheta =t_\hlam(\bX)$
when $L_{0,n}$ is replaced by the $L_{2,n}$ in (\ref{L_2n}).
\end{cor}

\def\theequation{5.\arabic{equation}}
\setcounter{equation}{0}
\section{Discussion}\label{discussion}

Although the focus of this paper is adaptive optimality sharp to the constant,
the oracle inequality in Theorem \ref{thm-4}
also implies the following rate optimality properties as corollaries.

\begin{cor}\label{cor-3} Let $\bX$, $\btheta$, $P_{\btheta}$, $t_\lam(x)$ and
$\hlam$ be as in Theorem \ref{thm-4}. Then,
\bel{cor-3-1}
 E_{\btheta}\|t_\hlam(\bX)-\btheta\|^2
\le M^*(1+\delta_{2,n}) \min_{\lam\ge 0}E_{\btheta}\|s_\lam(\bX)-\btheta\|^2
\eel
for all vectors $\btheta$ satisfying $\|\btheta\|_2^2\ge L_{2,n}/n^{1+\delta_{1,n}}$,
where $M^*$ is a constant depending on $\{\alpha_1',\alpha_2',\beta_0,A,M_0\}$ 
only and $L_{2,n}$ is as in (\ref{L_2n}).
\end{cor}

\begin{cor}\label{cor-4} Let $\bX$, $\btheta$, $P_{\btheta}$, $t_\lam(x)$  and
$\hlam$ be as in Theorem \ref{thm-4}, $\scrR(\Theta)$ be the minimax risk (\ref{minimax-risk}), and
\bes
 \Omega_n^{s,w} &=& \Big\{(p,C): 0 < p\le 2-c^{s,w},
C^{p} \ge (L_{2,n}/n^{\delta_{1,n}})^{p/2}n^{-1} \Big\},
\cr  \Omega_{0,n} &=& \Big\{(p,C): p=0,  C\ge 1/n \Big\},
\ees
with $c^{s,w}=c^{s}=0$ for strong balls and any $c^{s,w}=c^{w}\in (0,1)$ for weak balls.
Then, the adaptive rate minimaxity holds in the following sense:
\bes
\sup_{(p,C)\in\Omega_n^{s,w}\cup\Omega_{0,n}}
\f{\sup_{\btheta\in\Theta_{p,C,n}^{s,w}}E_{\btheta}\|t_{\hlambda}(\bX)-\btheta\|^2}
{\scrR(\Theta_{p,C,n}^{s,w})}\le M^*(1+\delta_{2,n}),
\ees
where $M^*$ and $L_{2,n}$ are as in Corollary \ref{cor-3}.
\end{cor}

For $0\le p\le 2$, the minimax rate in $\ell_p$ balls can be expressed as
\bel{minimax-rate}
\scrR(\Theta_{p,C,n}^{s,w})\asymp
\begin{cases} M_p^{s,w}\min\Big(nC^{p}\lam_{p,C,n}^{2-p}, n, (nC^{p})^{2/p}\Big), & 0<p\le 2
\cr nC\lam_{0,C,n}^{2}, & p=0<C
\cr 0,& p=0=C,
\end{cases}
\eel
where $\lam_{p,C,n} = 1\vee \sqrt{2\log(n\wedge(1/C^{p'})}$ with
$p' = p$ for $p>0$ and $p'=1$ for $p=0$. Here $a\asymp b$ means $a/b=O(1)$ and
this $O(1)$ is uniform in (\ref{minimax-rate}).
It follows from (\ref{minimax-0}) that for small $C^{p'}$, the constant factor
in (\ref{minimax-rate}) is accurate in the sense of its uniform validity when $\asymp$ is
replaced by $\approx$, provided that $p=0$ or $nC^{p}\lam_{p,C,n}^{2-p} \ll (nC^{p})^{2/p}$.
When $nC^{p}\lam_{p,C,n}^{2-p} \asymp (nC^{p})^{2/p}$, or equivalently
$nC^p \asymp \lam_{n,C,p}^p$, the constant factor in (\ref{minimax-rate}) is no longer accurate
for $p>0$ by (\ref{minimax-risk-formula}).
Moreover, for $p>0$, $nC^{p}\lam_{p,C,n}^{2-p}$ is of greater order than the minimax rate when
$(nC^{p})^{2/p} \ll nC^{p}\lam_{p,C,n}^{2-p}$.

In addition to the exact adaptive minimaxity literature discussed
earlier, adaptive rate minimaxity in $\ell_p$ balls was proved in
\cite{BirgeM01} for generalized $C_p$ when the minimax $\ell_2$ risk
is of no smaller order than $O(1)$, and in \cite{JohnstoneS04} for
EBThresh when the risk is of no smaller order than $(\log
n)^{2+(2-p)/2}$ and for a modified EBThresh when the risk is of no
smaller order than $\log n$, among many important contributions to
the problem. It follows from \cite{BirgeM01,Zhang05} that a hybrid
between the Fourier general empirical Bayes estimator and universal
soft threshold estimators is also adaptive rate minimax in $\ell_p$
balls when the minimax $\ell_2$ risk is of no smaller order than
$O(1)$.

The results in \cite{JohnstoneS04,AbramovichBDJ06} are valid for the
$\ell^q$ loss with $q\ge p$. It is unclear at the moment of this
writing if our analysis can be extended to hard threshold estimators
and the $\ell_q$ loss for $q>2$. The continuity of the soft
threshold estimator is a significant element in our analysis.

%

\vspace{10mm}

\appendix

\def\theequation{A.\arabic{equation}}
\setcounter{equation}{0}
\section*{Appendix}\label{proof}

\vspace{5mm}
\textsc{Proof of Lemma \ref{lm-Bayes-approx}.}
(i) The risk of soft thresholding $N(0,1)$,
\bes
R(0,\lam) = E(|N(0,1)|-\lam)_+^2 = 2\int_\lam^\infty (x-\lam)^2\varphi(x)dx,
\ees
is clearly decreasing in $\lam\in [0,\infty)$. To prove (\ref{R(0,lambda)-bounds}), we define
\bes
J_k(\lambda)=\int_0^\infty
u^k\exp(-u-u^2/(2\lambda^2))du.
\ees
With $u=\lam(x-\lam)$, we find that $\varphi(x) = \varphi(\lam)\exp(-u-u^2/(2\lambda^2))$ and
\bel{pf-lm-Bayes-1}
\lambda^3R(0,\lambda)/2
&=& \lambda^3\int_\lam^\infty (x-\lam)^2\varphi(x)dx = \varphi(\lam)J_2(\lam), 
\cr \lambda \Phi(-\lambda)
&=& \lambda \int_\lam^\infty \varphi(x)dx = \varphi(\lam)J_0(\lam).
\eel
Integrating by parts yields
\bel{pf-lm-Bayes-2}
(k+1)J_k(\lambda)&=&\int_0^\infty\exp(-x-x^2/(2\lambda^2))dx^{k+1}\cr
&=&\int_0^\infty(x^{k+1}+x^{k+2}/\lambda^2)\exp(-x-x^2/(2\lambda^2))dx\\ \nonumber
&=&J_{k+1}(\lambda)+J_{k+2}(\lambda)/\lambda^2.
\eel
It follows that $J_0(\lambda)=J_1(\lambda)+J_2(\lambda)/\lambda^2\ge
J_2(\lambda)/2+J_2(\lambda)/\lambda^2$, so that
\bes
J_2(\lambda)/J_0(\lambda)\le1/(1/2+1/\lambda^2).
\ees
In addition, (\ref{pf-lm-Bayes-2}) also implies that $J_3(\lambda)\le3J_2(\lambda)$, so that
\bes
\f{J_2(\lambda)}{J_0(\lambda)}
= \f{J_2(\lambda)}{J_2(\lambda)/2+J_3(\lambda)/(2\lambda^2)+J_2(\lambda)/\lambda^2}
\ge \f{1}{1/2+3/(2\lambda^2)+1/\lambda^2}.
\ees
We complete the proof of (\ref{R(0,lambda)-bounds}) by simple algebra after applying the above
two displayed inequalities to (\ref{pf-lm-Bayes-1}).

It follows from (\ref{risk-soft-threshold}) that for $\mu\neq 0$,
\bel{R(mu,lambda)}
R(\mu,\lambda) = E\Big(\mu^2I_{|Z+\mu|\le\lam}+(Z-\lam)^2I_{Z+\mu>\lam}+(Z+\lam)^2I_{Z-\mu<-\lam}\Big)
\eel
with $Z\sim N(0,1)$.
This implies $(\pa/\pa\mu)R(\mu,\lambda)=2\mu P\big\{|N(\mu,1)|\le\lambda\big\}$.
Since $R(\mu,\lam)$ is even in $\mu$ and $\lim_{\mu\to\infty}R(\mu,\lam)=\lam^2+1$, we find
\bes
R(\mu,\lam) \le R(0,\lam) + \mu^2\wedge (\lam^2+1),\quad \mu\wedge\lam\ge 0.
\ees
Integration of this inequality with $dG$ gives the first inequality in (\ref{R_G-upper-bound}).
Since $\rho_G(b)\le (b/a)^2\rho_G(a)$ for $0\le a\le b$ and ${B}_0\ge 4$ in (\ref{r_G}), for $\lam\ge 1$
the second inequality in (\ref{R_G-upper-bound}) follows from the first and (\ref{R(0,lambda)-bounds}) via
\bes
 \rho_G\big(\sqrt{\lambda^2+1}\big)+R(0,\lambda)
\le \big(1+\lam^{-2}\big)\rho_G(\lam)+\frac{4\Phi(-\lam)}{\lam^2}
\le \rho_G(\lam)+\frac{r_G(\lam)}{\lam^2}.
\ees
Let $g_0(\lam) = 4\Phi(-\lam)+\lam^2-1 - R(0,\lambda)$.
For $0\le \lam\le 1$, we have $1-\lam^2\ge (\lam^2+1)\wedge u^2 - 2(\lam^2\wedge u^2)$, so that
\bes
&& \rho_G(\lam)+r_G(\lam) - \Big\{\rho_G\big(\sqrt{\lambda^2+1}\big)+R(0,\lambda)\Big\}
\cr &= & \int\Big\{2(\lam^2\wedge u^2)+4\Phi(-\lam) - (\lambda^2+1)\wedge u^2-R(0,\lambda)\Big\}G(du)
\cr &\ge& g_0(\lam).
\ees
Since $g_0(0)=0$, $g_0'(0)=-4\varphi(0) + 2E|N(0,1)|=0$,
and $g_0''(\lam)=4\lam\varphi(\lam)+2 - 2\Phi(-\lam)\ge 0$,
we have $g_0(\lam)\ge 0$.
Thus, the second inequality in (\ref{R_G-upper-bound}) also holds for $0\le\lam\le 1$.

(ii) Since $\Phi(-\lam-\tau\lam) \le e^{-\lam^2(\tau+\tau^2/2)}\Phi(-\lam)$, (\ref{R(0,lambda)-bounds}) implies
\bel{pf-lm-Bayes-3}
 R(0,\lam)\ge \frac{4\Phi(-\lam)}{\lam^2+5}\ge {B}_0\Phi(-\lam-\tau_1\lam),
\eel
where $\tau_1$ is the solution of $\tau_1+\tau_1^2/2=\lam^{-2}\log\{(\lam^2+5){B}_0/4\}$.
For $\lam\ge 1$, this implies $\tau_1\le \lam^{-2}\{2\log\lam +\log(5{B}_0/4)\}$.
We also need a lower bound for the difference $R_G(\lam)-R(0,\lam)$.
Again, since $(\pa/\pa\mu)R(\mu,\lambda)=2\mu P\big\{|N(\mu,1)|\le\lambda\big\}$
by (\ref{R(mu,lambda)}), for $0\le M\le\lambda$ we have
\bes
&&\Big\{1-2\Phi(-M)\Big\}\Big\{\mu^2\wedge(\lambda-M)^2\Big\}\\
&=&\int_0^{|\mu|\wedge(\lambda-M)}2uP\Big\{\lambda-2M\le N(\lambda-M,1)\le\lambda\Big\}du\\
&\le&\int_0^{|\mu|\wedge(\lambda-M)}2uP\Big\{|N(\lambda-M,1)|\le\lambda\Big\}du\\
&\le&\int_0^{|\mu|\wedge(\lambda-M)}2uP\Big\{|N(u,1)|\le\lambda\Big\}du\\
&\le&R(\mu,\lambda)-R(0,\lambda).
\ees
Let $\Phi(-M)=1/(2\vee \lambda)$. Since $\Phi(-M)\le e^{-M^2/2}/2$, $M<\sqrt{2\log((\lambda/2)\vee 1)}$.
Integrating over $G(d{u})$, we find that
\bes
R_G(\lambda)-R(0,\lambda)
\ge \Big(1-\f{2}{\lambda\vee 2}\Big)\rho_G\big(\lambda-\sqrt{2\log((\lambda/2)\vee 1)}\big).
\ees
Let $\tau_2=\sqrt{2\log((\lambda/2)\vee 1)}/(\lam\vee e)$.
Since $\rho_G(a)\le\rho_G(b)\le(b/a)^2\rho_G(a)$ for $0<a\le b$, for $\lam\ge e$ we have
\bel{R_G-R(0,lambda)-expansion}
R_G(\lambda)-R(0,\lambda)
&\ge& \Big(1-\f{2}{\lambda\vee 2}\Big)\frac{(1-\tau_2)^2}{(1+\tau_1)^2}\rho_G(\lam+\lam\tau_1)
\cr &\ge& \Big(1+2\tau_2+M^*_0/\lam\Big)^{-1}\rho_G(\lam+\lam\tau_1)
\eel
with an $M^*_0$ depending on $B_0$ only.
We are allowed to incorporate higher order terms in $M^*_0/\lam$ for $\lam\ge e$
in (\ref{R_G-R(0,lambda)-expansion}) since
$\tau_1\le \lam^{-2}\log\{(\lam^2+5){B}_0/4\}$ and $\tau_2\le \lam^{-1}\sqrt{2\log \lam}$.
Combining (\ref{pf-lm-Bayes-3}) and (\ref{R_G-R(0,lambda)-expansion}), we find that for $\lam\ge e$,
\bes
\rho_G(1) &\le& \rho_G(\lam+\lam\tau_1) + {B}_0\Phi(-\lam-\lam\tau_1)
\cr &\le & \Big(1+ \lam^{-1}\sqrt{8\log \lam}+M^*_0/\lam\Big) R_G(\lambda).
\ees
For $\lambda_{G}\ge e$, this implies (\ref{Bayes-approx-Lambda}) with $\lam=\lam_G$
due to $\eta^*_G\le \rho_G(\lam+\lam\tau_1) + B_0\Phi(\lam+\lam\tau_1)$.
For $0\le \lambda_{G}\le e$, we have $\eta_G =R_G(\lambda_{G}) \ge R(0,e)>0$,
$\eta_G^* \le r_G(0)={B}_0/2$ and $\rho_G(1)\le 1$, so that (\ref{Bayes-approx-Lambda})
also holds. 
$\hfill\square$

\def\Btil{{\tilde B}}
\vspace{5mm}
\textsc{Proof of Lemma \ref{lm-rho_G(1)to0-part}.}
(i) We first proof the monotonicity of $\lam^*_G$ in $B_0$.
Let $\tlam^*_G=\argmin_\lam\{\rho_G+\Btil_0\Phi(-\lam)\}$ with $\Btil_0\le B_0$.
For all $\lam\le\tlam^*_G$,
\bes
\rho_G(\tlam_G^*) - \rho_G(\lam) \le \Btil_0\Big\{\Phi(-\lam)-\Phi(-\tlam^*_G)\Big\}
\le B_0\Big\{\Phi(-\lam)-\Phi(-\tlam^*_G)\Big\},
\ees
so that $\lam^*_G\ge\tlam_G^*$.
We assume without loss of generality that $B_0=4$ in the proof of (\ref{lm-rho_G(1)to0-part-1})
and (\ref{lm-R_G-r_G-bd}) since they only involve lower bounds for $\lam^*_G$.

Let $o(1)$ denote uniformly small quantity when $\eta_n$ is sufficiently small.
Suppose $\eta_G^*\le\eta_n$. By (\ref{r_G}), $4\Phi(-\lam^*_G)\le \eta_n$.
Since $(d/dt)\{(\sqrt{\pi/2}+t)\Phi(-t) - \varphi(t)\} = -\sqrt{\pi/2}\varphi(t) + \Phi(-t)\le 0$,
we have $\Phi(-t)\ge \varphi(t)/(\sqrt{\pi/2}+t)$. Thus, with $t=\lam_G^*$, we find
$\varphi(\lam^*_G)\le (\sqrt{\pi/2}+\lam^*_G)\eta_n/4$, or equivalently
\bes
(\lam_G^*)^2 + \log\Big(\Big(\sqrt{\pi/2}+\lam^*_G \Big)^2\Big)
\ge 2\log\big(1/\eta_n\big) +\log(8/\pi).
\ees
Since $(1-x)^{1/2} = 1/(1+x/2+O(x^2))$ for small $x$, this implies
\bes
\lam^*_G
&\ge & \Big(2\log\big(1/\eta_n\big)+\log(8/\pi) - \log(2\log(1/\eta_n))+o(1)\Big)^{1/2}
\cr & = & \sqrt{2\log(1/\eta_n)}\Big(1 + \frac{\log\log(1/\eta_n)-\log(4/\pi)+o(1)}{4\log(1/\eta_n)}\Big)^{-1}.
\ees
This and (\ref{eta_G-bd}) implies $R(0,\lam_G)\le \eta_G\le (1+o(1))\eta_n$.
It follows from (\ref{R(0,lambda)-bounds}) that
$R(0,t) \ge 4\Phi(-t)/(t^2+5)\ge 4\varphi(t)/\{(t^2+5)(t+\sqrt{\pi/2})\}$,
so that
\bes
\lam_G
&\ge & \Big(2\log\Big(\frac{4+o(1)}{\eta_n\sqrt{2\pi}}\Big) -
2\log\Big\{\Big(\lam_G^2+5\Big)\Big(\lam_G +\sqrt{\pi/2}\Big)\Big\}\Big)^{1/2}
\cr & = & \Big(2\log(1/\eta_n)+\log(8/\pi)+o(1) - \log\lam_G^6\Big)^{1/2}
\cr &\ge & \sqrt{2\log(1/\eta_n)}\Big(1 + \frac{3\log\log(1/\eta_n)+\log \pi+o(1)}{4\log(1/\eta_n)}\Big)^{-1}.
\ees
This and (\ref{Bayes-approx-Lambda}) implies $\eta^*_G\le (1+o(1))\eta_n$.
This completes the proof of (\ref{lm-rho_G(1)to0-part-1}).

Consequently, the first inequality in (\ref{lm-R_G-r_G-bd}) follows from
(\ref{eta_G-bd}), (\ref{lm-rho_G(1)to0-part-1}) and the fact that $\rho_G(\lam)\le r_G(\lam)$,
while the second inequality in (\ref{lm-R_G-r_G-bd}) follows from
(\ref{Bayes-approx-Lambda}) and (\ref{lm-rho_G(1)to0-part-1}).

(ii) Let $\lam_n=\sqrt{2\log(1/\eta_n)}$.
Since $\lam_n\ge 4/\sqrt{2\pi}\ge 1$, $\rho_G(\lam_n)\le\lam_n^2\eta_n$ and
$\rho_G(\sqrt{\lam_n^2+1})\le (1+\lam_n^{-2})\rho_G(\lam_n)\le \rho_G(\lam_n)+\rho_G(1)$.
We also have
\bes
{B}_0\Phi(-\lam_n)\le ({B}_0/\lam_n)\varphi(\lam_n)=\eta_n {B}_0/(\lam_n\sqrt{2\pi})\le\eta_n.
\ees
Thus, by (\ref{r_G-bd}), $R_G(\lam_n) - \rho_G(1)\le \rho_G(\lam_n)+{B}_0\Phi(-\lam_n) \le (1+\lam_n^2)\eta_n$.


(iii) Let $t_n = 2\log\big(4n/\sqrt{2\pi}\big)$ and $y_n = t_n-3\log t_n$.
For $n\ge 10$, $t_n$ is increasing in $n$ with $3\log t_n > 5.13$ and
$y_n$ is increasing in $t_n$ with $y_n\{\exp(2/(3y_n))-1\}\le (2/3)\exp(2/(3y_n))\le 3.5\le 3\log t_n$. By algebra,
\bes
y_n + 3\log y_n + 2/y_n < t_n= y_n + 3\log t_n,\quad \forall \, n\ge 10.
\ees
By the Jensen inequality, $\int_0^\infty e^{-u-t u^2/2}du > e^{-t}$,
so that
\bes
\frac{1}{4n} = z_n^{-2}\Phi( -z_n) = z_n^{-3}\varphi(z_n)\int_0^\infty e^{-u-(u/z_n)^2/2}du > z_n^{-3}\varphi(z_n) e^{-1/z_n^2}.
\ees
This gives $z_n^2+3\log z_n^2 +2/z_n^2 > t_n$, so that $z_n^2 > y_n$ for $n\ge 10$.
We also numerically verify $z_n^2 > y_n$ for $7\le n\le 9$.

For $n\ge 1500$, $(\pa/\pa n)(y_n-\log n) = 2/n - (3/t_n)(2/n) - 1/n = (1-6/t_n)/n > 0.0004>0$
and $y_n-\log n \ge 0.01$, so that $z_n^2\ge y_n> \log n$. We also verify $z_n^2> \log n$
numerically for $2\le n < 1500$.

Suppose $\eta^*_{G_n}\le z_n^2/n$. Since $z_n^{-2}\Phi(-z_n) = 1/(4n)$, we have
\bes
B_0\Phi(-\lam^*_{G_n}) \le \eta^*_{G_n} \le 4\Phi\Big(-\sqrt{n \eta^*_{G_n}}\Big).
\ees
This inequality and the constraint $B_0\ge 4$ yield
\bes
\lam_{G_n}^*\ge \sqrt{n \eta^*_{G_n}}
> \sqrt{n\rho_{G_n}(\lam^*_{G_n})}
= \sqrt{\sum_{i=1}^n(\theta_i)^2\wedge(\lam^*_{G_n})^2}.
\ees
Consequently, $\rho_{G_n}(\lam^*_{G_n})=\|\btheta\|^2/n = r_{G_n}(\infty) = \eta_{G_n}^*$
and $\lam^*_{G_n}=\infty$.

Let $\eta_n=(\sqrt{2\log n}-2\sqrt{2\log\log n})^2/n$ with sufficiently large $n\ge n^*$.
Suppose $\eta_{G_n}\le \eta_n$. We need to prove $\eta^*_{G_n}\le z_n^2/n$.
Since $\lam^*_{G_n}$ is increasing in $B_0$, it suffices to consider $B_0=4$.
Since $\log(1/\eta_n) \ge \log n - \log(2\log n)$
it follows from (\ref{lm-R_G-r_G-bd}) and the condition $\eta_{G_n}\le\eta_n$ that
\bes
n \eta_{G_n}^*
&\le& \left(1+ \frac{2\sqrt{\log\log n}+M^*}{\sqrt{\log n}}\right)n\eta_n
\cr & = & \left(1+ \frac{2\sqrt{\log\log n}+M^*}{\sqrt{\log n}}\right)
\left(1- \frac{2\sqrt{\log\log n}}{\sqrt{\log n}}\right)^2 2\log n
\cr & = & \left(1- \frac{(2+o(1))\sqrt{\log\log n}}{\sqrt{\log n}}\right)2\log n
\cr & = & 2\log n - (4+o(1))\sqrt{\log\log n}\sqrt{\log n}.
\ees
Since $z_n^2\ge 2\log\big(4n/\sqrt{2\pi}\big) - 3\log\big(2\log\big(4n/\sqrt{2\pi}\big)\big)
= 2\log n - 3\log\log n +O(1)$, we have $n\eta^*_{G_n}\le z_n^2$.
$\hfill\square$

\vspace{5mm}
\textsc{Proof of Lemma \ref{lm-firm}.}  Assume $\mu\ge 0$ by symmetry.

(i) Let $x_{\mu,\lam}$ be the solution of $f_\lam(x_{\mu,\lam})=2\mu$.
Since $f_\lam(x) \le \kappa_0(x-\lam)$ for $x>0$, we have $(2-\kappa_0)\mu\le \kappa_0(x_{\mu,\lam}-\mu-\lam)$.
For $t_\lam(x)>2\mu$,
\bes
t_\lam(x)-\mu &\le& \kappa_0(x-\lam)-\mu
\cr &=& \kappa_0(x-\mu-\lam)+(\kappa_0-1)\mu
\cr &\le & \kappa_0(x-\mu-\lam)+\frac{(\kappa_0-1)}{2-\kappa_0}\kappa_0(x_{\mu,\lam}-\mu-\lam).
\ees
Since $f_\lam(x)\ge t_\lam(x)>2\mu$, we have $x>x_{\mu,\lam}$, so that
\bes
|t_\lam(x)-\mu| \le C_0(x-\mu-\lam) \le C_0(|x-\mu|-\lam)_+.
\ees
For $x>-\lam$ and $t_\lam(x)\le 2\mu$, $|t_\lam(x)-\mu|\le\mu$.
For $x \le -\lam$,
\bes
&& |t_\lam(x)-\mu|
\le \mu - f_\lam(x) \le \mu + \kappa_0(|x|-\lam) \le \kappa_0(|x-\mu|-\lam).
\ees
Thus, $|t_\lam(x)-\mu|$ is bounded by either $|\mu|$ or $C_0(|x-\mu|-\lam)_+$.

(ii) It suffices to consider the location model where the distribution of $X-\mu$ is fixed.
Since $(\pa/\pa\mu) E\big(s_\lam(X)-\mu\big)^2=2\mu P\{|X|\le\lam\}$,
$E\big(s_\lam(X)-\mu\big)^2\le \lim_{\mu\to\infty}E\big(s_\lam(X)-\mu\big)^2
=E(X-\mu-\lam)^2=\lam^2+\Var(X)$.
Since the value of $t_\lam(x)$ is between $s_\lam(x)$ and $x$, it follows that
\bes
E\Big(t_\lam(X)-\mu\Big)^2 \le E\Big(s_\lam(X)-\mu\Big)^2+E(X-\mu)^2
\le \lam^2+2\,\Var(X).
\ees

(iii) The $\ell_q$ error bound follows from (i) and (ii).

(iv) Since $B_0\ge 2C_0^2$, (\ref{R(0,lambda)-bounds}) gives
$C_0^2E(|X_i-\theta_i|-\lam)_+^2\le C_0^2 R(0,\lam)\le 4C_0^2 \Phi(-\lam)/\lam^2 \le 2B_0\Phi(-\lam)/\lam^2$.
Thus, (iii) with $q=2$ and (\ref{compound}) yield
\bes
n^{-1} E_{\btheta}\|t_\lam(\bX)-\btheta\|_2^2
&\le& \int \min\Big(u^2+C_0^2 R(0,\lam),\lam^2+2\Big)G_n(du)
\cr &\le& (1+2/\lam^2)\rho_{G_n}(\lam) + C_0^2 R(0,\lam)
\cr &\le& \rho_{G_n}(\lam) + 2r_{G_n}(\lam)/\lam^2.
\ees
For $0\le\lam\le 1$, $(\lam^2+2)\wedge u^2 - 3(\lam^2\wedge u^2)\le 2(1-\lam^2)\le (B_0/2)(1-\lam^2)$ gives
\bes
&& \rho_{G_n}(\lam) + 2r_{G_n}(\lam) - \rho_{G_n}(\sqrt{\lam^2+2}) - C_0^2R(0,-\lam)
\cr &\ge & \int\Big\{3(\lam^2\wedge u^2)+2B_0\Phi(-\lam) - (\lam^2+2)\wedge u^2 - (B_0/2)R(0,\lam)\Big\}G_n(du)
\cr &\ge& 2B_0\Phi(-\lam) + (B_0/2)(\lam^2 - 1) - (B_0/2)R(0,\lam),
\ees
which is nonnegative as in the proof of (\ref{r_G-bd}).
$\hfill\square$

\vspace{5mm}
\textsc{Proof of Lemma \ref{lm-Gauss-mix}.}
(i) Let $h_G(t)=\int e^{-u^2/2}\cosh(t{u})G(d{u})$ where
$\cosh(t)=(e^t+e^{-t})/2$. Since
$P\big\{|N(u,1)|>t\big\}=1-\Phi(t-u)+\Phi(-t-u)$,
\bes
\f{S_G(t)}{2\Phi(-t)}&=&\f{1}{2\Phi(-t)}
\int\int_{x>t}\Big(\varphi(x+u)+\varphi(x-u)\Big)dxG(du)\\
&=&\f{1}{2\Phi(-t)}\int\int_{x>t}\Big(e^{-u^2/2}(e^{-xu}+e^{xu})\Big)\varphi(x)dxG(du)\\
&=&\f{1}{\Phi(-t)}\int_t^\infty h_G(x)\varphi(x)dx.
\ees
Since $G$ does not put the entire mass at $0$,
 $h_G(t)$ is strictly increasing in $t$ for $t\ge 0$, and
the monotonicity of $S_G(t)/\Phi(-t)$ follows from 
\bes &&\f{\pa}{\pa
t}\log\Big(\f{S_G(t)}{2\Phi(-t)}\Big)\\
&=&-\f{h_G(t)\varphi(t)}{\int_t^\infty h_G(x)\varphi(x)dx}+
\f{\varphi(t)}{\Phi(-t)}\\
&=&\f{\varphi(t)}{\Phi(-t)\int_t^\infty
h_G(x)\varphi(x)dx}\Big(\int_t^\infty
h_G(x)\varphi(x)dx-h_G(t)\Phi(-t)\Big) > 0.
\ees

Since
$P\big\{|N(\mu,1)|>t\big\}$ is even in $\mu$ and \bes
\Big(\f{\pa}{\pa\mu}\Big)^2P\big\{|N(\mu,1)|>t\big\}
&=&-\Phi''(t-\mu)+\Phi''(-t-\mu)\\
&\le&2\max_tt\varphi(t)=2\varphi(1)<2, \ees we have
$S_G(t)\le\int\big(P\big\{|N(0,1)|>t\big\}+1\wedge
u^2\big)G(d{u})=2\Phi(-t)+\rho_G(1)$. In addition, \bes
\Gbar(t)/2\le\int_{|u|\ge t}P\big\{N(|u|,1)>t\big\}G(du)\le
S_G(t).\ees These imply the inequalities in (\ref{Gauss-mix-1}).

To prove $\lam^*_{G_n} > \xi_{2,*}$, we observe from (\ref{r_G}) that
\bel{G(lambda*)}
\frac{ dr_{G_n}(\lam)}{d\lam}\Big|_{\lam=\lam^*_{G_n}}
= 2\lam^*_{G_n}\Gbar_n(\lam^*_{G_n}) - {B}_0\varphi(\lam^*_{G_n})=0,
\eel
so that by (\ref{Gauss-mix-1}) and the condition ${B}_0\ge 8/\alpha_2'$
\bes
\f{S_{G_n}(\lam^*_{G_n})}{2\Phi(-\lam^*_{G_n})}
\ge \f{\Gbar_n(\lam^*_{G_n})/2}{2\Phi(-\lam^*_{G_n})}
= \f{{B}_0\varphi(\lam^*_{G_n})}{8\lam^*_{G_n}\Phi(-\lam^*_{G_n})} > \frac{1}{\alpha_2'}.
\ees
The monotonicity of $S_{G_n}(t)/\Phi(-t)$ guarantees $\lam^*_{G_n}>\xi_{2,*}$ by (\ref{xi_*}).

(ii) Since $S_{G_n}(\xi_{j,*})/(2\Phi(-\xi_{j,*}))=1/\alpha'_j$,
$1/\alpha_j'-1\le \rho_{G_n}(1)/(2\Phi(-\xi_{j,*}))$ by (\ref{Gauss-mix-1}) and simply algebra.
For $t\le\xi_{j,*}$, we have  $g_{1,n}(t)\le t\le \xi_{j,*}$ and $R(0,g_{1,n}(t))\le 4\Phi(-t)$ by (\ref{g_1-cond}).
Thus, by the monotonicity of $S_{G_n}(t)/\Phi(-t)$, (\ref{xi_*}) and (\ref{R(0,lambda)-bounds}),
(\ref{Gauss-mix-2}) follows from
\bes
&& \frac{\alpha_j'S_{G_n}(g_{1,n}(t))}{2\Phi(-t)}\le \frac{\Phi(-g_{1,n}(t))}{\Phi(-t)}
\le \frac{(5+g_{1,n}^2(t))R(0,g_{1,n}(t))}{4\Phi(-t)} \le 5+t^2.
\ees

(iii) Since $\theta_{*,n}^2 \le 2\beta_0\log n \le\log n$, $\theta_{*,n} \le z_n$ in
Lemma \ref{lm-rho_G(1)to0-part} (iii), so that $\|\btheta\|^2 = n \eta_{G_n}^* \le \theta_{*,n}^2$.
It follows that $\|\btheta\|_\infty\le\theta_{*,n}$ and $\int|u| G_n(du)\le \sqrt{\|\btheta\|_2^2/n} \le \theta_{*,n}/\sqrt{n}$
by Cauchy-Schwarz.
By the convexity of $\Phi(x)$ in $x<0$,
\bel{pf-10}
S_{G_n}(t) &=&\ \int\big\{\Phi(-t-|u|)+\Phi(-t+|u|)\big\}G_n(du)
\cr &\le & \Phi(-t)+\int\bigg\{\Big(1-\f{|u|}{\theta_{*,n}}\Big)\Phi(-t) +\f{|u|}{\theta_{*,n}}\Phi(-t+\theta_{*,n})\bigg\}G_n(du)
\\ \nonumber & \le & 2\Phi(-t)+n^{-1/2}\Big\{\Phi(-t+\theta_{*,n})-\Phi(-t)\Big\}
\eel
for $t\ge \theta_{*,n}$. For $\theta_{*,n}\le t\le A_1^{1/2}\xi_{1,1}$, $(t\theta_{*,n}+\theta_{*,n}^2/2)\le\beta_0\log n$, so that
\bes
\Phi(-t+\theta_{*,n})-\Phi(-t)
\le \Phi(-t)e^{t\theta_{*,n}+\theta_{*,n}^2/2}
\le \Phi(-t)n^{\beta_0}.
\ees
Consequently, $S_{G_n}(t)/\{2\Phi(-t)\} \le 1+n^{\beta_0-1/2}/2 \le 1/\alpha_1'$ at $t=A_1^{1/2}\xi_{1,1}$.
This gives $\xi_{1,*}\ge A_1^{1/2}\xi_{1,1}$ by (\ref{xi_*}) and the monotonicity of $S_{G_n}(t)/\Phi(-t)$.

The proof of (\ref{Gauss-mix-3}) utilizes the following fact.
For $x>t>0$, $\Phi(-x) =\int_t^\infty \varphi(u+x-t)du\le e^{-t(x-t)-(x-t)^2/2}\Phi(-t)$, so that
\bes
\int_t^\infty \Phi^{1/2}(-x)dx
&\le& \int_t^\infty \Phi^{1/2}(-t)e^{-t(x-t)/2-(t-x)^2/4}dx
\cr &\le& \Phi^{1/2}(-t)\min\big(2/t,\sqrt{\pi}\big).
\ees
Since $t-\theta_{*,n}\ge t/2$ for $2\theta_{*,n}\le t\le A_1^{1/2}\xi_{1,1}$, (\ref{pf-10}) implies
\bes
\int_t^\infty S_{G_n}^{1/2}(x)dx
&\le& \int_t^\infty \sqrt{2\Phi(-x) + n^{-1/2}\Phi(-x+\theta_{*,n})}dx
\cr &\le& \sqrt{2\Phi(-t)}2/t + \sqrt{n^{-1/2}\Phi(-t+\theta_{*,n})}2/(t-\theta_{*,n})
\cr &\le& \sqrt{2\Phi(-t)}(2/t)\Big(1 + 2\sqrt{n^{-1/2}+n^{\beta_0-1/2}}\Big).
\ees
This completes the proof of the lemma.
$\hfill\square$

\vspace{5mm}
\textsc{Proof of Lemma \ref{lm-exp-inq}.}
(i) Let $t_k=\log(k/E_{\btheta}N(\xi_k))$.
For $k>E_{\btheta}N(\xi_k)$, we have $t_k>0$ and
\bes
P_{\btheta}\big\{N(\xi_k)\ge k\big\}&\le&e^{-t_kk}E_{\btheta}\exp(t_kN(\xi_k))\\
&=&e^{-t_kk}\prod_{i=1}^nE_{\btheta}\exp(t_kI\big\{|X_i|\ge\xi_k\big\})\\
&=&e^{-t_kk}\prod_{i=1}^n\Big(1+(e^{t_k}-1)P_{\btheta}\big\{|X_i|\ge
\xi_k\big\}\Big)\\
&\le&\exp\{-t_kk+(e^{t_k}-1)E_{\btheta}N(\xi_k)\}\\
&=&\exp(-\nu_kk),
\ees
where $\nu_k=E_{\btheta}N(\xi_k)/k-1-\log(E_{\btheta}N(\xi_k)/k)>0$.
Similarly, for $k<E_{\btheta}N(\xi_k)$ and $t_k<0$, \bes
P_{\btheta}\big\{N(\xi_k)\le
k\big\}&\le&e^{-t_kk}E_{\btheta}\exp(t_kN(\xi_k))\\
&\le&\exp\{-t_kk+(e^{t_k}-1)E_{\btheta}N(\xi_k)\}=\exp(-\nu_kk).\ees
Thus, (\ref{exp-inq-1}) holds in both cases.

(ii) Due to the monotonicity of $S_{G_n}(t)/\Phi(-t)$, for $\xi_{1,k}\le \xi_{1,*}$
\bes
&&\f{S_{G_n}(\xi_{1,k})}{k/n}=\f{\alpha_1S_{G_n}(\xi_{1,k})}{2\Phi(-\xi_{1,k})}
\le \f{\alpha_1S_{G_n}(\xi_{1,*})}{2\Phi(-\xi_{1,*})} = \f{\alpha_1}{\alpha'_1}.
\ees
Since $\alpha_1<\alpha_1'<1$, we have $k/n>S_{G_n}(\xi_{1,k})$, so that by (\ref{exp-inq-1})
\bes
P_{\btheta}\big\{N(\xi_{1,k})-k\ge0\big\}\le\exp(-\nu_{1,k}k).
\ees
Since $x-1-\log(x)$ is a decreasing function for $0<x<1$,
and $E_{\btheta}N(\lambda)/n=S_{G_n}(\lambda)$ by the definition these quantities in
(\ref{N(t)}) and (\ref{S_G-Gbar}),
\bes
\nu_{1,k} &=&\f{S_{G_n}(\xi_{1,k})}{k/n}-1-\log\Big(\f{S_{G_n}(\xi_{1,k})}{k/n}\Big)\\
&\ge&\f{\alpha_1}{\alpha'_1}-1-\log\Big(\f{\alpha_1}{\alpha'_1}\Big)=\nu_{1,*}.
\ees
The above inequalities imply (\ref{exp-inq-2}) in view of the definition of $\hxi_1$ in (\ref{xi}).
The proof of (\ref{exp-inq-3}) is nearly identical and omitted.
$\hfill\square$

\vspace{5mm}
\textsc{Proof of Lemma \ref{lm-Gaussian-iso-inq2}.}
Let $N(t)$ be as in (\ref{N(t)}) and define
\bes
\Delta(\bx;a,b)=\max_{a \le \lambda\le b}\|t_\lambda(\bx)-t_{b}(\bx)\|/\sqrt{n}.
\ees
The following inequalities follow directly from related definitions and (\ref{smooth-cond}):
\bel{Delta-inq}
&&\sup_{a\le\lambda\le b}\Big|L(\bX,\btheta,\lambda)-L(\bX,\btheta,b)\Big|\le\Delta(\bX;a,b),\cr 
&&\sup_{a\le\lambda\le b}\Big|R_{1,n}(\btheta,\lambda)-R_{1,n}(\btheta,b)\Big|
\le E_{\btheta}\Delta(\bX;a,b),\\ \nonumber
&& E_{\btheta}\Delta^2(\bx;a,b) \le \kappa_1^2(b-a)^2E_{\btheta}N(a)/n = \kappa_1^2(b-a)^2S_{G_n}(a).
\eel
The last inequality in (\ref{Delta-inq}) follows from
$|t_\lambda(x)-t_{b}(x)|\le\kappa_1(b-a)I\big\{|x|>a\big\}$ for $a\le\lambda\le b$.

Let $\Delta_j(\bx)=\Delta(\bx;c_{j-1},c_j)$ and $B_j=\big\{c_{j-1}\le Y < c_j\big\}$.
It follows from the triangle inequality and the first two inequalities of (\ref{Delta-inq}) that
\bel{pf-lm-concen-inq-discretization}
\sqrt{H_{\btheta}(a,b)}
&\le&\max_{1\le j\le m}\Big|L(\bX,\btheta,Y)-R_{1,n}(\btheta,Y)\Big|I_{B_j}\nonumber\\
&\le&\max_{1\le j\le m}\Big|L(\bX,\btheta,c_j)-R_{1,n}(\btheta,c_j)\Big|I_{B_j}\nonumber\\
&&+\max_{1\le j\le m}\Big\{\Delta_j(\bX)+E_{\btheta}\Delta_j(\bX)\Big\} I_{B_j} \\
&\le&\max_{1\le j\le m}\Big|L(\bX,\btheta,c_j)-R_{1,n}(\btheta,c_j)\Big|I_{B_j}\nonumber\\
&&+\max_{1\le j\le m}\Big\{\Delta_j(\bX)-E_{\btheta}\Delta_j(\bX)\Big\}_+ I_{B_j}
 +2\max_{1\le j\le m}I_{B_j}E_{\btheta}\Delta_j(\bX). \nonumber
\eel
Since $\|t_\lam(\bX)-\btheta\|$ also has Lipschitz norm $\kappa_0$, $L(\bx,\btheta,\lam)$ has
Lipschitz norm $\kappa_0/\sqrt{n}$. Thus, by the Gaussian isoperimetric inequality,
\bes
&&E_{\btheta}\bigg\{\max_{1\le j\le
m}\Big|L(\bX,\btheta,c_j)-R_{1,n}(\btheta,c_j)\Big|I_{B_j}\bigg\}^2\\
&\le&\int_0^\infty\sum_{j=1}^mP_{\btheta}
\Big\{\big|L(\bX,\btheta,c_j)-R_{1,n}(\btheta,c_j)\big|I_{B_j}>x\Big\}dx^2\\
&\le& \sum_{j=1}^m\int_0^\infty \min\big(\pi_j, 2e^{-nx^2/(2\kappa_0^2)}\big)dx^2\\
&=&\f{2\kappa_0^2}{n}\sum_{j=1}^m\bigg\{ \pi_j \log(2/\pi_j)+\pi_j \bigg\}.
\ees
Similarly, due to $|\Delta_j(\bu)-\Delta_j(\bv)|\le\kappa_0\|\bu-\bv\|/\sqrt{n}$,
\bes
&&E_{\btheta}\bigg\{\max_{1\le
j\le m}\Big(\Delta_j(\bX)-E_{\btheta}\Delta_j(\bX)\Big)_
+I_{B_j}\bigg\}^2\le\f{2\kappa_0^2}{n}\sum_{j=1}^m\pi_j\log(2e/\pi_j).
\ees
Inserting the above two inequalities and the third inequality of (\ref{Delta-inq}) to
(\ref{pf-lm-concen-inq-discretization}) after an application of the Minkowski inequality,
we find
\bes
\sqrt{E_{\btheta} H_{\btheta}(a,b)}
&\le& \bigg\{E_{\btheta}\Big(\max_{1\le j\le m}
\Big|L(\bX,\btheta,c_j)-R_{1,n}(\btheta,c_j)\Big|I_{B_j}\Big)^2\bigg\}^{1/2}
\cr &&+\bigg\{E_{\btheta}\Big(\max_{1\le j\le
m}\Big(\Delta_j(\bX)-E_{\btheta}\Delta_j(\bX)\Big)_+
I_{B_j}\Big)^2\bigg\}^{1/2}
\cr &&+2\bigg\{E_{\btheta}\Big(\max_{1\le j\le
m}I_{B_j}E_{\btheta}\Delta_j(\bX)\Big)^2\bigg\}^{1/2}
\cr &\le&2\kappa_0\bigg\{\f{2}{n}\sum_{j=1}^m\pi_j\log(2e/\pi_j)\bigg\}^{1/2}
 +2\kappa_1\bigg\{\sum_{j=1}^m\pi_j (c_j-c_{j-1})^2S_{G_n}(c_{j-1})\bigg\}^{1/2}.
\ees
This is (\ref{lm-Gaussian-iso-discrete}).
With $c_j=a+(j/m)(b-a)$ and $\pi_j=P\big\{c_{j-1}\le Y <c_j\big\}$,
\bes
\sum_{j=1}^m \pi_j\log(2e/\pi_j)
&=& \sum_{j=1}^m \int_0^\infty \min\Big\{\pi_j, 2e^{-t}\Big\}dt
\cr &\le& \int_0^\infty\min\Big\{1,2m e^{-t}\Big\}dt
= \log(2em)
\ees
and $(c_j-c_{j-1})^2S_{G_n}(c_{j-1})\le m^{-2}(b-a)^2S_{G_n}(a)$,
so that (\ref{lm-Gaussian-iso-discrete}) implies (\ref{lm-Gaussian-iso-discrete-coll}).

It follows from the definition of $H^*_{\btheta}(a,b)$ and (\ref{Delta-inq}) that for $\epsilon>0$,
\bes
\sqrt{E_{\btheta}H^*_{\btheta}(\lam)} - \sqrt{E_{\btheta}H^*_{\btheta}(\lam+\epsilon)}
\le \sqrt{E_{\btheta}\Delta^2(\bX;\lam,\lam+\epsilon)}\le \epsilon \kappa_1S_{G_n}^{1/2}(\lam),
\ees
in view of the third inequality of (\ref{Delta-inq}).
This implies $(\pa/\pa \lam) \sqrt{E_{\btheta}H^*_{\btheta}(\lam)} \le \kappa_1S_{G_n}^{1/2}(\lam)$.
Since $H^*_{\btheta}(\lam)\to 0$ almost surely as $\lam\to\infty$, the monotone convergence
theorem gives (\ref{lm-Gaussian-iso-star}).

Finally, the Gaussian isoperimetric inequality gives \bes
\Var[L(\bX,\btheta,\lambda)]&=&\int_0^\infty
P_{\btheta}\Big\{\big|L(\bX,\btheta,\lambda)-R_{1,n}(\btheta,\lambda)\big|>t\Big\}dt^2\\
&\le&\f{1}{n}\int_0^\infty
P_{\btheta}\Big\{\Big|\|t_\lambda(\bX)-\btheta\|-E_{\btheta}\|t_\lambda(\bX)-\btheta\|\Big|>t\Big\}dt^2\\
&\le&\f{2}{n}\int_0^\infty e^{-t^2/(2\kappa_0^2)}dt^2=\f{4\kappa_0^2}{n}.\ees This and
$E_{\btheta}L^2(\bX,\btheta,\lambda)=R_{1,n}^2(\btheta,\lambda)+\Var[L(\bX,\btheta,\lambda)]$
give (\ref{lm-Gaussian-iso-basic}).
$\hfill\square$

\vspace{5mm}
\textsc{Proof of (\ref{cond-g_1-a}).} Since $R(0,x)$ is decreasing in $t$, (\ref{cond-g_1}) implies
$g_{1,n}(x)\to\infty$ as $x\to\infty$. For large $x$, $R(0,g_{1,n}(x))=(4+o(1))g_{1,n}(x)^{-3}\varphi(g_{1,n}(x))
\le 4\Phi(-x)=(4+o(1))x^{-1}\varphi(x)$ by (\ref{cond-g_1}), so that
\bes
- g_{1,n}^2(x)/2 - 3\log g_{1,n}(x) +o(1) \le - x^2/2 - \log x.
\ees
This implies $g_{1,n}(x)\ge (1+o(1))x$. It follows that  $g_{1,n}(x) = (1+o(1))x$ due to  $g_{1,n}(x)\le x$.
Thus, for $A_1=1+\delta_{1,n}$ and large $x$
\bes
&& \Phi\Big(-\sqrt{1+\delta_{1,n}}g_{1,n}(x)\Big)
\cr &\approx & (\sqrt{A_12\pi}x)^{-1}\left(\sqrt{2\pi}\varphi(g_{1,n}(x))\right)^{1+\delta_{1,n}}
\cr &\approx & \frac{(2\pi)^{\delta_{1,n}/2}}{A_1^{1/2}x}\left((x^3/4)R(0,g_{1,n}(x))\right)^{1+\delta_{1,n}}
\cr &\le & \frac{(1+o(1))(2\pi)^{\delta_{1,n}/2}}{x}
\left(\frac{x^3M_0\Phi(-x)}{4(x^{c_{1,n}}+2)(\log_+x)^{c_{2,n}}}\right)^{1+\delta_{1,n}}.
\ees
This completes the proof of (\ref{cond-g_1-a}). $\hfill\square$

\vspace{5mm}
\textsc{Proof of Theorem \ref{thm-4}.}
Let $A_j = 1+\delta_{j,n}$, $j=1,2$. We denote by $M^*$ a constant depending on
$\{\alpha_1',\alpha_2',\beta_0,A,M_0\}$ only which may take different values
from one appearance to the next.
We note that $A_1\le A$ for all $n$ and $A_2 \le A$ for $n< n_*$.

Recall that in (\ref{lam-star}), $\lam_{G_n}^*=\argmin_\lam r_{G_n}(\lam)$ and
$\eta_{G_n}^*=r_{G_n}(\lam_{G_n}^*)$. Our plan is to prove that
\bel{pf-thm-4-1}
&& E_{\btheta}\Big(\|t_\hlam(\bX)-\btheta\|/\sqrt{n} - R_{1,n}(\btheta,\hlam)\Big)^2
\le M^*\Big(A_2\tau^*_{1,n}\eta_{G_n}^*+\tau^*_{2,n}\Big)
\eel
and that with $R^{(sm)}_{G_n}(\lam)=E_{\btheta}\|t_\lam(\bX)-\btheta\|^2/n$ as in (\ref{R-sm}),
\bel{pf-thm-4-2}
 E_{\btheta}R^{(sm)}_{G_n}(\hlam)
\le A_2\eta_{G_n}^*
+ M^*\Big(A_2\tau^*_{1,n}\eta_{G_n}^*+\tau^*_{2,n}\Big).
\eel
We first observe that (\ref{thm-4-1}) follows from (\ref{pf-thm-4-1}) and (\ref{pf-thm-4-2}); To wit,
\bes
&& \sqrt{E_{\btheta}\|t_\hlam(\bX)-\btheta\|^2/n}
\cr &\le& \sqrt{E_{\btheta} R_{1,n}^2(\btheta,\hlam)}
+\sqrt{E_{\btheta}\Big(\|t_\hlam(\bX)-\btheta\|/\sqrt{n} - R_{1,n}(\btheta,\hlam)\Big)^2}
\cr &\le& \sqrt{A_2\eta_{G_n}^*
+ M^*\big(A_2\tau^*_{1,n}\eta_{G_n}^*+\tau^*_{2,n}\big)}
+\sqrt{M^*(A_2\tau^*_{1,n}\eta_{G_n}^*+\tau^*_{2,n})}
\cr &\le& \sqrt{A_2\eta_{G_n}^*}
+2\sqrt{M^*(A_2\tau^*_{1,n}\eta_{G_n}^*+\tau^*_{2,n})},
\ees
due to the Cauchy-Schwarz inequality
$R_{1,n}^2(\btheta,\lam) \le R^{(sm)}_{G_n}(\lam)$.

Let 
$\theta_{*,n} =\max\{t: A_1^{1/2}\xi_{1,1}t+t^2/2\le \beta_0(1\vee\log n)\}$.
Define
\bel{cases}
\begin{cases}
 \hbox{ Case 1:}\ \ n < n_*\ \hbox{ or }\ \eta_{G_n}^* > \theta_{*,n}^2/n,
\cr \hbox{ Case 2:}\ \ n\ge n_*\ \hbox{ and }\ \eta_{G_n}^*\le \theta_{*,n}^2/n.
\end{cases}
\eel
It follows from (\ref{xi_jk}) and the definition of $\theta_{*,n}$ that
\bes
\frac{1}{\theta_{*,n}^2}
= \frac{\big(A_1^{1/2}\xi_{1,1} + \sqrt{A_1\xi_{1,1}^2 + 2\beta_0\log n}\big)^2}{(2\beta_0)^2(1\vee \log n)^2}
\le \frac{M^*A_1}{1\vee \log n}.
\ees
Thus, due to $
\log(e\vee \log n)/(1\vee \log n) \le \tau^*_{1,n}$ by (\ref{tau_n}) and
the boundedness of $\max_{n<n_*}1/\tau_{2,n}^*=\max_{n<n_*} n^{A_1}/L_{2,n}$ by (\ref{L_2n}),
we have in Case 1
\bel{Case-1-bd}
\frac{1}{n}
\le \begin{cases}
\eta_{G_n}^*/\theta_{*,n}^2 \le
M^*\tau^*_{1,n}\eta_{G_n}^*/\log(e\vee \log n), & n\ge n_*
\cr M^*\tau^*_{2,n}/\{A_2n\log(e\vee \log n)\}, & n<n_*.
\end{cases}
\eel
The conditions in Case 2 allows application of 
Lemma \ref{lm-Gauss-mix} (iii). 

Let $\nu_{1,*}=\alpha_1/\alpha'_1-1-\log(\alpha_1/\alpha'_1)$ as in (\ref{exp-inq-2}) and define
\bel{k_0}
 k_0 = \begin{cases} 1,& \hbox{\ Case 1}, \cr
\big(\lceil(\delta_{1,n}\log n+2\log\log n)/\nu_{1,*}\rceil \vee 1\big) \wedge n,& \hbox{\ Case 2}.
\end{cases}
\eel
Let $\xi_{1,*}$ be as in (\ref{xi_*}). Define $\xi_{1,n+1} = \xi_{1,*}$ for $\xi_{1,*}<\xi_{1,n}$,
\bel{lambda_*}
k_{1,*} &=& \min\big\{k: \xi_{1,k} \le (\xi_{1,*}\wedge\xi_{1,k_0})\big\},
\cr \lambda_{1,*} &=& \sqrt{1+\delta_{1,n}}g_{1,n}(\xi_{1,k_{1,*}})=A_1^{1/2}g_{1,n}(\xi_{1,k_{1,*}}).
\eel
Note that $k_{1,*}\ge k_0$ and that $\{k_{1,*},\xi_{1,k_{1,*}}\}$ is the unique solution of
\bes
\begin{cases}
\xi_{1,k_{1,*}} \le \xi_{1,*} < \xi_{1,k_{1,*}-1}, & \xi_{1,n}\le \xi_{1,*} \le \xi_{1,k_0}\cr
\xi_{1,k_{1,*}}=\xi_{1,k_0}, & \xi_{1,*} \ge \xi_{1,k_0} \cr
\xi_{1,k_{1,*}}=\xi_{1,*}, & \xi_{1,*} < \xi_{1,n}.
\end{cases}
\ees

We split the excess risk in 4 main terms:
\bel{zeta_1n}
&\zeta_{1,n} \equiv E_{\btheta}\left|\|t_{\hlambda}(\bX)-\btheta\|/\sqrt{n}
-R_{1,n}(\btheta,\hlambda)\right|^2I\big\{\hlambda\le \lambda_{1,*}\big\},
\\ \label{zeta_2n}
&\zeta_{2,n} \equiv E_{\btheta}\left|\|t_{\hlambda}(\bX)-\btheta\|/\sqrt{n}
-R_{1,n}(\btheta,\hlambda)\right|^2I\big\{\hlambda > \lambda_{1,*}\big\},
\\ \label{zeta_3n}
&\zeta_{3,n} \equiv E_{\btheta}\,\rho_{G_n}((\hlam^2+2)^{1/2}) - \rho_{G_n}(\lam_{G_n}^*),
\\ \label{zeta_4n}
&\zeta_{4,n} \equiv C_0^2E_{\btheta}\,R(0,\hlam).
\eel
We prove in four steps that
\bel{zeta_1n-bd}
&\zeta_{1,n} \le M^*\Big(\tau^*_{1,n}\eta_{G_n}^*+\tau^*_{2,n}\Big),
\\ \label{zeta_2n-bd}
& \displaystyle \zeta_{2,n} \le M^*\Big(A_2\tau^*_{1,n}\eta_{G_n}^*+\tau^*_{2,n}\Big),
\\ \label{zeta_3n-bd}
&\zeta_{3,n} \le \delta_{2,n}\eta_{G_n}^* + A_2 M^*\tau^*_{1,n}\eta_{G_n}^*,
\\ \label{zeta_4n-bd}
&\zeta_{4,n} \le M^*\Big(\tau^*_{1,n}\eta_{G_n}^*+\tau^*_{2,n}\Big).
\eel
Inequalities (\ref{zeta_1n-bd}) and (\ref{zeta_2n-bd}) directly imply (\ref{pf-thm-4-1}).
By (\ref{r_G}) and (\ref{lam-star}), $\eta_{G_n}^*=\rho_{G_n}(\lam_{G_n}^*)+{B}_0\Phi(-\lam_{G_n}^*)$.
By Lemma \ref{lm-firm} (iv),
\bes
R^{(sm)}_G(\hlam) \le \rho_G((\hlam^2+2)^{1/2}) + C_0^2R(0,\hlam).
\ees
Thus, (\ref{zeta_3n-bd}) and (\ref{zeta_4n-bd}) imply (\ref{pf-thm-4-2}).
It remains to prove (\ref{zeta_1n-bd}), (\ref{zeta_2n-bd}), (\ref{zeta_3n-bd}) and (\ref{zeta_4n-bd}).
This is done in the following four steps respectively.

\vspace{5mm}
\textsc{Step 1.}
In this step we prove (\ref{zeta_1n-bd}).
Since $g'_{1,n}(x) > 0$, $\xi_{1,*}< \xi_{1,n}$ implies
$\hlam \ge A_1^{1/2}g_{1,n}(\xi_{1,n}) > A_1^{1/2}g_{1,n}(\xi_{1,*})=\lam_{1,*}$.
Thus, this step only concerns the case of $\xi_{1,*}\ge \xi_{1,n}$,
where $\xi_{1,n}\le \xi_{1,k_{1,*}}\le \xi_{1,*}$.

It follows from (\ref{hlam}) and (\ref{exp-inq-2}) that for all $k \ge k_{1,*}$,
\bes 
 P\Big\{g_{1,n}(\xi_{1,k+1})\le A_1^{-1/2}\hlambda \le g_{1,n}(\xi_{1,k})\Big\}
\le P\Big\{\hxi_1 \le \xi_{1,k}\Big\}
\le e^{-\nu_{1,*}k}.
\ees
Since $2\Phi(-\xi_{j,k}) =\alpha_jk/n$, we have
\bel{pf-th1-step1-2a}
\alpha_j/(2n) &=& \Phi(-\xi_{j,k+1})-\Phi(-\xi_{j,k})
\cr &\ge& \varphi(\xi_{j,k})(\xi_{j,k}-\xi_{j,k+1})
\cr &\ge&  \xi_{j,k}\Phi(-\xi_{j,k})(\xi_{j,k}-\xi_{j,k+1})
\\ \nonumber &=&  \xi_{j,k}(\xi_{j,k}-\xi_{j,k+1})\alpha_jk/(2n).
\eel
Since $0\le (d/dx) g_{1,n}(x)\le M_0$ by (\ref{cond-g_1}), this gives
\bes
&& 0 \le g_{1,n}(\xi_{1,k})-g_{1,n}(\xi_{1,k+1})\le M_0(\xi_{1,k}-\xi_{j,k+1}) \le M_0/(k\xi_{1,k}).
\ees
An application of (\ref{lm-Gaussian-iso-discrete}) with $\pi_k=e^{-\nu_{1,*}k}$ and $c_{n+1-k}=A_1^{1/2}g_{1,n}(\xi_{1,k})$ yields
\bel{pf-th1-step1-2}
&&E_{\btheta}\bigg|\f{1}{\sqrt{n}}\|\hbtheta-\btheta\|-R_{1,n}(\btheta,\hlambda)\bigg|^2
I\big\{\hlambda \le \lambda_{1,*}\big\}\\ \nonumber
&\le& 8\sum_{k_{1,*}\le k < n}\bigg\{\f{\pi_k\log(2e/\pi_k)}{n/(2\kappa_0^2)}
+\pi_k\Big(\frac{A_1^{1/2}M_0}{k\xi_{1,k}/\kappa_1}\Big)^2S_{G_n}\big(g_{1,n}(\xi_{1,k+1})\big)\bigg\}.
\eel
For $k_{1,*}\le k < n$, we have $\xi_{1,k+1}\le\xi_{1,*}$ by (\ref{lambda_*}),
so that by (\ref{Gauss-mix-2})
\bes
\frac{S_{G_n}(g_{1,n}(\xi_{1,k+1}))}{\alpha_1(k+1)/n}
= \frac{S_{G_n}(g_{1,n}(\xi_{1,k+1}))}{2\Phi(-\xi_{1,k+1})}
\le \f{5+\xi_{1,k+1}^2}{\alpha_1'}\le \f{5+\xi_{1,k}^2}{\alpha_1'}.
\ees
It follows from this inequality, (\ref{zeta_1n}) and (\ref{pf-th1-step1-2}) that
\bel{step-1-fin}
\zeta_{1,n} \le \sum_{k_{1,*}\le k}e^{-\nu_{1,*}k}
\bigg\{\f{k\nu_{1,*}+\log(2e)}{n/(16\kappa_0^2)}
+\f{A_1M_0^2(5+\xi_{1,k}^2)(k+1)}{k^2\xi_{1,k}^2n \alpha_1'/(8\alpha_1\kappa_1^2)}\bigg\}
\le \frac{M^*k_{1,*}\nu_{1,*}}{ne^{\nu_{1,*}k_{1,*}}}.
\eel
We note here that $1/\xi_{1,k}\le 1/\xi_{1,n} = 1/|\Phi^{-1}(\alpha_1/2)|$ is bounded.

In Case 1, (\ref{zeta_1n-bd}) follows from (\ref{step-1-fin}) and (\ref{Case-1-bd})
due to $x/e^{x}\le 1/e$:
\bes
\zeta_{1,n} \le \frac{M^*}{n} \le M^*\Big(\tau^*_{1,n}\eta^*_{G_n}+\tau^*_{2,n}/n\Big).
\ees

In Case 2, $\nu_{1,*}k_{1,*}\ge \nu_{1,*}k_0 \ge \delta_{1,n}\log n + 2\log\log n$
by (\ref{k_0}) and (\ref{lambda_*}), so that (\ref{step-1-fin}) implies
\bes
\zeta_{1,n}\le \frac{M^*k_0}{n^{1+\delta_{1,n}}(\log n)^2}\le \frac{M^*}{n^{1+\delta_{1,n}}\log n}
\le M^*\tau^*_{2,n}
\ees
by the definition of the $\tau^*_{2,n}$. Note that by (\ref{L_2n})
\bes
&& n^{1+\delta_{1,n}}\tau^*_{2,n}= L_{2,n}\ge (\log\log n)^{1-c_{2,n}}/(\log n)^{c_{1,n}/2}\ge \frac{\log\log n}{M^*\log n},
\ees
due to the constraints $c_{2,n}\le 0$ for $c_{1,n}=2$ and $0 < c_{1,n}\le 2$.

\vspace{5mm}
\textsc{Step 2.}
In this step, we prove (\ref{zeta_2n-bd}).
We first consider Case 1 as specified in (\ref{cases}) with $k_0=1$ as in (\ref{k_0}).
An application of the concentration inequality (\ref{lm-Gaussian-iso-discrete-coll}) yields
that for all positive integer $m$,
\bes
\zeta_{2,n}&=&E_{\btheta}\bigg|\f{1}{\sqrt{n}}\|\hbtheta-\btheta\|
-R_{1,n}(\btheta,\hlambda)\bigg|^2I\Big\{\lambda_{1,*}\le\hlambda\le \sqrt{A_2}\xi_{2,1}\Big\}\\
&\le& M^*\Big(\f{\log m}{n}+S_{G_n}(\lambda_{1,*})\f{A_2\xi_{2,1}^2}{m^2}\Big).
\ees
By (\ref{xi_jk}), $\xi_{1,k_{1,*}}^2\le \xi_{1,1}^2\le \xi_{2,1}^2\le 2\log(2n/\alpha_2)$.
By the upper bound for $S_{G_n}(g_{1,n}(t))$ in (\ref{Gauss-mix-2}),
\bes
S_{G_n}(\lambda_{1,*}) \le (5+\xi_{1,k_{1,*}}^2)2\Phi(-\xi_{1,k_{1,*}})/\alpha_1'
\le M^*(\log n)2\Phi(-\xi_{1,k_{1,*}}).
\ees
Since $\xi_{1,k_{1,*}}\le \xi_{1,*}\le \xi_{1,k_{1,*}-1} \le \xi_{1,k_0}$ for $k_{1,*}>k_0=1$
and $2\Phi(-\xi_{1,k})=\alpha k/n$,
the upper bound for $2\Phi(-\xi_{j,*})$ in Lemma \ref{lm-Gauss-mix} (ii) yields
\bes
2\Phi(-\xi_{1,k_{1,*}})
\le \frac{\alpha_1}{n}+4\Phi(-\xi_{1,*})
\le \frac{\alpha_1}{n} + \frac{2\alpha_1'\rho_{G_n}(1)}{1-\alpha_1' }.
\ees
Moreover, (\ref{Bayes-approx-Lambda}) 
and (\ref{lm-R_G-r_G-bd}) imply $\rho_{G_n}(1)\le M^*\eta_{G_n}^*$.
Thus, with $m=\lceil (\log n)^2\rceil$ and an application of (\ref{Case-1-bd}),
the upper bound for $\zeta_{2,n}$ becomes
\bes
\zeta_{2,n} &\le& M^*\Big(\f{\log m}{n}+S_{G_n}(\lambda_{1,*})\f{A_2\xi_{2,1}^2}{m^2}\Big)\\
&\le& M^*\Big\{\f{\log\log n}{n} +\f{A_2(\log n)^2}{(\log n)^4}
\Big(\frac{1}{n}+\rho_{G_n}(1)\Big)\Big\}\\
&\le& M^*\Big\{\f{A_2+\log\log n}{n} +\f{A_2\eta_{G_n}^*}{(\log n)^2}\Big\}\\
&\le& M^*\big(A_2\tau^*_{1,n}\eta_{G_n}^*+\tau^*_{2,n}\big).
\ees

\vspace{3mm}
Now we consider Case 2 as specified in (\ref{cases}).
By Lemma \ref{lm-Gauss-mix} (iii),
$\xi_{1,*}\ge A_1^{1/2}\xi_{1,1}\ge \xi_{1,k_0}$, so that $\lam_{1,*} = A_1^{1/2} g_{1,n}(\xi_{1,k_0})$ by (\ref{lambda_*}).
Since $2\Phi(- \xi_{1,k_0}) = \alpha_1k_0/n$ and $\log k_0 = o(1)\log n$ by (\ref{k_0}),
$\xi_{1,k_0}\approx \sqrt{2\log n}$.
By (\ref{cond-g_1}), $(1+o(1))x = g_{1,n}(x)\le x$ for large $x$ as in the proof of (\ref{cond-g_1-a}),
so that $\lam_{1,*}\approx A_1\sqrt{2\log n}$.
Thus, by (\ref{cond-g_1-a}) and the definition of $L_{2,n}$ in (\ref{L_2n}),
\bes
\Phi(-\lam_{1,*}) &= & \Phi\Big(-\sqrt{1+\delta_{1,n}}g_{1,n}(\xi_{1,k_0})\Big)
\cr &\le & \frac{M^*}{\sqrt{1+\delta_{1,n}}\xi_{1,k_0}}
\left(\frac{\xi_{1,k_0}^3\Phi(-\xi_{1,k_0})}{\xi_{1,k_0}^{c_{1,n}}(\log\xi_{1,k_0})^{c_{2,n}}}\right)^{1+\delta_{1,n}}
\cr &\le & \frac{M^*}{\sqrt{1+\delta_{1,n}}\xi_{1,1}}
\left(\frac{\xi_{1,1}^3\alpha_1k_0/n}{\xi_{1,1}^{c_{1,n}}(\log\xi_{1,1})^{c_{2,n}}}\right)^{1+\delta_{1,n}}
\cr &\le & \frac{M^*}{\sqrt{\log n}}
\left(\frac{(\log n)^{(3-c_{1,n})/2}k_0}{n(\log\log n)^{c_{2,n}}}\right)^{1+\delta_{1,n}}
\cr &\le &  \frac{M^*}{\sqrt{\log n}}
\left(\frac{\delta_{1,n}(\log n)^{(5-c_{1,n})/2}}{n(\log\log n)^{c_{2,n}}}
+\frac{(\log n)^{(3-c_{1,n})/2}}{n(\log\log n)^{c_{2,n}-1}} \right)^{1+\delta_{1,n}}
\cr &\le & \frac{M^*L_{2,n}}{n^{1+\delta_{1,n}}} \lam_{1,*}^2.
\ees
Moreover, it follows from the definition of $H^*_{\btheta}(\lam)$ in (\ref{H-star}) that
\bes
&& \left|\|t_{\hlambda}(\bX)-\btheta\|/\sqrt{n}
-R_{1,n}(\btheta,\hlambda)\right| I\big\{\hlambda > \lambda_{1,*}\big\}
\cr &\le& \sqrt{H^*_{\btheta}(\lambda_{1,*})} + E_{\btheta}\sqrt{ H^*_{\btheta}(\lambda_{1,*})}.
\ees
Consequently, (\ref{lm-Gaussian-iso-star}) of Lemma \ref {lm-Gaussian-iso-inq2}
and (\ref{Gauss-mix-3}) of Lemma \ref{lm-Gauss-mix} (iii) yield
\bel{step-2-fin}
\frac{\zeta_{2,n}}{\kappa_1} \le 4\Big(\int_{\lam_{1,*}}^\infty S_{G_n}^{1/2}(t)dt\Big)^2
\le \frac{(4C_j')^2}{\lam_{1,*}^2}\Phi(-\lam_{1,*}) \le \frac{M^*L_{2,n}}{n^{1+\delta_{1,n}}}=M^*\tau^*_{2,n}.
\eel
Thus, (\ref{zeta_2n-bd}) holds in both cases.

\vspace{5mm}
\textsc{Step 3.}
In this step we prove (\ref{zeta_3n-bd}). Recall that $\rho_G(\lam)=\int(u^2\wedge \lam^2)G(du)$ in (\ref{r_G})
and $\hlam\le \sqrt{1+\delta_{2,n}}\hxi_2$ in (\ref{hlam}).
We bound $\zeta_{3,n}$ in (\ref{zeta_3n}) by
\bel{step3-1}
\zeta_{3,n} \le \Gbar_n(\lam^*_{G_n})
E_{\btheta}\Big(2+(1+\delta_{2,n})\hxi_2^2-(\lam^*_{G_n})^2\Big)_+.
\eel
By (\ref{pf-th1-step1-2a}), $\xi_{2,k}^2-\xi_{2,k+1}^2\le (\xi_{2,k}+\xi_{2,k+1})/(k\xi_{2,k})\le 2/k$.
Since ${B}_0 = 8/\alpha_2'$, $\lam^*_{G_n}\ge \xi_{2,*}$ by Lemma \ref{lm-Gauss-mix} (i).
Let $k_{2,*}=\inf\{k\ge 1:\xi_{2,k} \ge \lam^*_{G_n}\}$.
By (\ref{exp-inq-3}),
\bes
E_{\btheta}\Big(\hxi_2^2-(\lam^*_{G_n})^2\Big)_+
&\le & E_{\btheta}\Big(\hxi_2 - \xi_{2,k_{2,*}+1}^2\Big)_+
\cr &=& \sum_{1 \le k \le k_{2,*}}(\xi_{2,k}^2-\xi_{2,k+1}^2) P_{\btheta}\big\{\hxi_2 \ge \xi_{2,k}\big\}
\cr & \le & \sum_{1 \le k \le k_{2,*}}\frac{2}{k} e^{-\nu_{2,*}k}.
\ees
In view of (\ref{step3-1}) and the definitions in (\ref{r_G}) and (\ref{lam-star}), it follows that
\bes
\zeta_{3,n} &\le & \Gbar_n(\lam^*_{G_n})\Big\{2+(1+\delta_{2,n})M^*
+\delta_{2,n}(\lam^*_{G_n})^2\Big\}
\cr & \le & \eta_{G_n}^*\Big\{\delta_{2,n}+(1+\delta_{2,n})M^*/(\lam^*_{G_n})^2\Big\}.
\ees
It follows from (\ref{G(lambda*)}) that
$\lam^*_{G_n}\ge \lam^*_{G_n}\Gbar_n(\lam^*_{G_n}) = (B_0/2)\varphi(\lam^*_{G_n})$,
so that $1/\lam_{G_n}^*$ is uniformly bounded.
Thus, $1/(\lam^*_{G_n})^2 \le M^*/\log_+(1/\eta_{G_n}^*)\le M^*\tau_{1,n}^*$
by (\ref{lm-rho_G(1)to0-part-1}) of Lemma \ref{lm-rho_G(1)to0-part} (i) and the definition of $\tau_{1,n}^*$,
and (\ref{zeta_3n-bd}) follows.

\vspace{5mm}
\textsc{Step 4.} In this step we prove (\ref{zeta_4n-bd}).
By Lemma \ref{lm-Bayes-approx} (i), $R(0,\lam)$ is decreasing in $\lam$.
Let $\hlam_1 = A_1^{1/2}g_{1,n}(\hxi_1)$. Since $\hlam\ge \hlam_1$, it suffices to prove
\bel{step-4-0}
E_{\btheta}R\left(0,\hlam_1\right) \le M^*\big(\tau^*_{1,n}\eta_{G_n}^*+\tau^*_{2,n}\big).
\eel

By condition (\ref{cond-g_1}) on $g_{1,n}$ and the definition of $\xi_{1,k}$ in (\ref{xi_jk}),
$R(0,g_{1,n}(\xi_{1,k}))\le 4\Phi(-\xi_{1,k}) = 2\alpha_1k/n$.
By (\ref{lambda_*}) and (\ref{exp-inq-2}) of Lemma \ref{lm-exp-inq},  
\bes
E_{\btheta}R\left(0,\hlam_1\right) I\big\{\hlam_1\le\lam_{1,*}\big\}
&\le& \sum_{k_{1,*}\le k \le n} R\left(0,g_{1,n}(\xi_{1,k})\right) P_{\btheta}\big\{\hxi_1=\xi_{1,k}\big\}
\\ \nonumber &\le& \sum_{k_{1,*}\le k \le n} 2\alpha_1 (k/n) e^{-\nu_{1,*}k}
\ees
Thus, by the analysis in Step 1 after (\ref{step-1-fin}),
\bel{step-4-1}
&& E_{\btheta}R\left(0,\hlam_1\right) I\big\{\hlam_1\le\lam_{1,*}\big\}
\le M^*(\tau^*_{1,n}\eta_{G_n}^*+\tau^*_{2,n}).
\eel

Since $R(0,\lam)$ is decreasing in $\lam$ and $\lam_{1,*} = \sqrt{1+\delta_{1,n}}g_{1,n}(\xi_{1,k_{1,*}})$,
\bel{step4-2}
E_{\btheta}R\left(0,\hlam_1\right) I\big\{\hlam_1 > \lam_{1,*}\big\}
&\le & R\left(0,\lam_{1,*}\right)
\\ \nonumber &\le& \begin{cases} R(0,g_{1,n}(\xi_{1,*}\wedge \xi_{1,1})), & \hbox{Case 1} \cr
R(0,\sqrt{1+\delta_{1,n}}g_{1,n}(\xi_{1,k_0})), & \hbox{Case 2}
\end{cases}
\eel
with the two cases specified in (\ref{cases}) and the definition of $k_{1,*}$ in (\ref{lambda_*}).
Note that we proved $\xi_{1,*}\ge \xi_{1,1}\ge \xi_{1,k_0}$ in Case 2 in Step 2.

We first consider Case 1.
If $\eta_{G_n}^*\ge 1/M^*$, then $\tau^*_{1,n}$ is not small and
$R\left(0,g_{1,n}(\xi_{1,*})\right)\le 1\le M^*\eta_{G_n}^*\le (M^*)^2\tau^*_{1,n}\eta_{G_n}^*$.
Otherwise, $\eta_{G_n}^*\le 1/M^*$, so that the upper bound for $2\Phi(-\xi_{j,*})$
in Lemma \ref{lm-Gauss-mix} (ii), (\ref{lm-rho_G(1)to0-part-1}) and (\ref{lm-R_G-r_G-bd}) imply
\bes
2\Phi(-\xi_{1,*}) \le \frac{\alpha_1'\rho_{G_n}(1)}{1-\alpha_1'} \le M^* \eta_{G_n}^*.
\ees
Thus, as in the proof of (\ref{lm-rho_G(1)to0-part-1}),
\bes
 \sqrt{2\log(1/\eta_{G_n}^*)}
 \le \left(1+ \frac{\log\log(1/\eta_{G_n}^*)+M^*}{4\log(1/\eta_{G_n}^*)}\right)\xi_{1,*}\le 2\xi_{1,*}.
\ees
The second bound for $R(0,g_{1,n}(x))$ in (\ref{cond-g_1}) gives
\bes
R\left(0,g_{1,n}(\xi_{1,*})\right)
&\le& \frac{M_0\Phi(-\xi_{1,*})}{(\xi_{1,*}^{c_{1,n}}+2)(\log_+ \xi_{1,*})^{c_{2,n}}}
\cr &\le& \frac{M^*\eta_{G_n}^*}{(\log_+(1/\eta_{G_n}^*))^{c_{1,n}/2}(\log_+\log(1/\eta_{G_n}^*))^{c_{2,n}}}
\cr &\le& M^*\tau^*_{1,n}\eta_{G_n}^*,
\ees
in view of the definition of $\tau^*_{1,n}$ in (\ref{tau_n}).
Moreover, the first bound for  $R(0,g_{1,n}(x))$ in (\ref{cond-g_1}) and (\ref{Case-1-bd}) give
\bes
R\left(0,g_{1,n}(\xi_{1,1})\right)\le 4\Phi(-\xi_{1,1})=\frac{2\alpha_1}{n}\le M^*(\tau^*_{1,n}\eta_{G_n}^*+\tau^*_{2,n}).
\ees
It follows that $R\left(0,g_{1,n}(\xi_{1,*}\wedge \xi_{1,1})\right)\le M^*(\tau^*_{1,n}\eta_{G_n}^*+\tau^*_{2,n})$,
and (\ref{zeta_4n-bd}) follows from (\ref{step4-2}).

In Case 2, (\ref{R(0,lambda)-bounds}) of Lemma \ref{lm-Bayes-approx} (i) and (\ref{step-2-fin}) yield
\bes
R\left(0,\lam_{1,*}\right)
\le \frac{4\Phi(-\lam_{1,*})}{\lam_{1,*}^2+2}
\le M^*\tau^*_{2,n}.
\ees
Thus, 
(\ref{step-4-0}) holds in both cases in view of (\ref{step-4-1}) and (\ref{step4-2}).
The proof of Theorem \ref{thm-4} is completed
since we have already proved the oracle inequality based on
(\ref{zeta_1n-bd}), (\ref{zeta_2n-bd}), (\ref{zeta_3n-bd}) and (\ref{zeta_4n-bd})
$\hfill\square$


\vspace{5mm}
\textsc{Proof of Theorems \ref{th-adaptive-ratio-opt} and \ref{th-smooth-thresh} and Corollary \ref{cor-1}.}
It follows from (\ref{lm-R_G-r_G-bd}) that $\eta_{G_n}^*\le (1+M^*\tau_{1,n})\eta_{G_n}$,
so that (\ref{cor-1-1}) follows from Theorem~\ref{thm-4}.
Since $M_nL_{2,n}/n^{1+\delta_{1,n}}\le\eta_n\to 0$, we have $n\to\infty$
and $\tau_{1,n}\to 0$.
The adaptive ratio optimality (\ref{def-adaptive-ratio-opt}) then follows from (\ref{cor-1-1}) with
the special $t_\lam(x)=s_\lam(x)$
since the risk range guarantees
$\tau^*_{2,n} \ll \eta_{G_n}$ uniformly in the specified class.
Theorems \ref{th-adaptive-ratio-opt} and \ref{th-smooth-thresh} are consequences of Corollary \ref{cor-1}
since (\ref{cond-g_1}) holds and $L_{2,n}=L_{0,n}$ for $g_{1,n}(x)=x$, $c_{1,n}=2$ and $c_{2,n}=0$.
$\hfill\square$

\vspace{5mm}
\textsc{Proof of Theorem \ref{th-adaptive-minimax} and Corollary \ref{cor-2}.}
Let $\Omega_n^{s,w}$ and $\Omega_{0,n}$ be as in Theorem~\ref{th-adaptive-minimax}
with $L_{0,n}$ replaced by $L_{2,n}$.
It follows from the second part of (\ref{minimax-risk-formula}),
which implies (\ref{minimax-0}), that for certain $M_n\to\infty$ and $\eta_n\to 0$,
\bes
M_n L_{2,n}/n^{\delta_{1,n}} \le \scrR(\Theta_{p,C,n}^{s,w}) \le n\eta_n
\ees
uniformly for all $(p,C)\in\Omega_n^{s,w}\cup\Omega_{0,n}$.
Thus, Corollary \ref{cor-1} and the first part of (\ref{minimax-risk-formula}) imply that uniformly
for all $(p,C)\in\Omega_n^{s,w}\cup\Omega_{0,n}$,
\bes
&& \sup\Big\{E_{\btheta}\|t_\hlam(\bX)-\btheta\|^2: \btheta\in\Theta_{p,C,n}^{s.w}\Big\}
\cr &\le & \sup\Big\{n\Big(\sqrt{\eta_{G_n}} + M^*\sqrt{\tau_{1,n}\eta_{G_n} + \tau^*_{2,n}}\Big)^2:
\btheta\in\Theta_{p,C,n}^{s.w}\Big\}
\cr &\le & (1+o(1))\sup\Big\{n\eta_{G_n}:
\btheta\in\Theta_{p,C,n}^{s.w}\Big\}
\cr & = & (1+o(1))\scrR(\Theta_{p,C,n}^{s,w}).
\ees
This completes the proof of Corollary \ref{cor-2}.
Theorem \ref{th-adaptive-minimax} is a consequence of Corollary \ref{cor-2}
since (\ref{cond-g_1}) holds for $g_{1,n}(x)=x$, $c_{1,n}=2$ and $c_{2,n}=0$
and $L_{2,n}=L_{0,n}$ for those $c_{1,n}$ and $c_{2,n}$.
$\hfill\square$

\vspace{5mm}
The proof of Theorem \ref{th-min-loss} requires the following lemma.

\begin{lem}\label{lm-bound-threshold-diff}
For any real numbers $\lambda>b\ge0$, $\mu$ and $\eps$, \bes
\big(s_b(\eps+\mu)-\mu\big)^2-\big(s_{\lambda}(\eps+\mu)-\mu\big)^2
\le\big(|\eps|+b\big)^2I\big\{|\eps|>b\big\}.\ees
\end{lem}

\textsc{Proof.} Let $\mu>0$ without loss of generality due to
symmetry. We have \bes
&&\big(s_b(\eps+\mu)-\mu\big)^2-\big(s_{\lambda}(\eps+\mu)-\mu\big)^2\\
&\le&\left\{
\begin{array}{ll}
(\eps-b)^2 &\mbox{if $\eps>b$},\\
(\eps+b)^2 &\mbox{if $\eps+\mu\le-b$,}\\
0 &\mbox{if $|\eps+\mu|\le b$,}\\
(\eps-b)^2-\mu^2\le0 &\mbox{if $b<\eps+\mu\le\lambda$ and $\eps\le b$,}\\
(\eps-b)^2-(\eps-\lambda)^2\le0 &\mbox{if $\eps+\mu>\lambda$ and
$\eps\le b$.}
\end{array}
\right.\ees The upper bound is no greater than
$\big(|\eps|+b\big)^2I\big\{|\eps|>b\big\}$. $\hfill\square$

\vspace{5mm}
\textsc{Proof of Theorem \ref{th-min-loss}.}
Let $a=0< b=\sqrt{2\log n}$ and
\bes
Y_1=\argmin_{\lambda\ge0}\|s_{\lambda}(\bX)-\btheta\|^2,\quad
Y_2=\min(Y_1,b).
\ees
The risk difference between $s_{Y_1}$ and $s_{Y_2}$ is
controlled by applying Lemma \ref{lm-bound-threshold-diff},
\bel{bound-Y1-Y2}
&&\f{1}{n}E_{\btheta}\|s_{Y_2}(\bX)-\btheta\|^2-\f{1}{n}E_{\btheta}\|s_{Y_1}(\bX)-\btheta\|^2
\cr &\le& \f{1}{n}\sum_{i=1}^nE_{\btheta}\Big\{\big(|X_i-\theta_i|+b\big)^2I\big\{|X_i-\theta_i|>b\big\}\Big\}
\cr &\le& 2\int_b^\infty(t+b)^2\varphi(t)dt
\cr &\le& M^*\Big(\f{\log n}{n}\Big).
\eel
Since $0\le Y_2\le b=\sqrt{2\log n}$, an application of the concentration inequality
(\ref{lm-Gaussian-iso-discrete-coll}) with $m=n$ gives that
\bes
&& \sqrt{E_{\btheta}\big|\|s_{Y_2}(\bX)-\btheta\|/\sqrt{n}-R_{1,n}(\btheta,Y_2)\big|^2}
\cr &\le& 2\bigg\{\f{2}{n}\big(\log(2en)+1\big)\bigg\}^{1/2} +2\bigg\{\f{2\log n}{n^2}\bigg\}^{1/2}
\\ \nonumber &\le& M^*\sqrt{(\log n)/n}.
\ees
This and (\ref{bound-Y1-Y2}) yield
\bel{pf-thm3-Schwarz}
&& \sqrt{E_{\btheta}R_{1,n}^2(\btheta,Y_2)}
\le M^*\sqrt{(\log n)/n} + \sqrt{E_{\btheta}\|s_{Y_1}(\bX)-\btheta\|^2/n}.
\eel
It follows from the optimality of $\lam_{G_n}= \argmin_\lam R_{G_n}(\lam)$
and (\ref{lm-Gaussian-iso-basic}) that
\bes
R_{G_n}(\lam_{G_n}) \le E_{\btheta}R_{G_n}(Y_2) \le E_{\btheta}R_{1,n}^2(\btheta,Y_2) + 4/n.
\ees
Consequently,
\bes
\sqrt{R_{G_n}(\lam_{G_n})}
&\le& \sqrt{4/n} + M^*\sqrt{(\log n)/n} + \sqrt{E_{\btheta}\|s_{Y_1}(\bX)-\btheta\|^2/n}
\cr &=& \sqrt{4/n} + M^*\sqrt{(\log n)/n} + \sqrt{E_{\btheta}\inf_{\lam}\|s_{\lam}(\bX)-\btheta\|^2/n}.
\ees
This and (\ref{adaptive-ratio-optimality}) yield Theorem \ref{th-min-loss}. $\hfill\square$


\vspace{1.5cm}

\bibliographystyle{abbrv}
\bibliography{FDR-thresh}

\bigskip
\noindent
Wenhua Jiang\\
School of Mathematical Sciences\\
Soochow University\\
P.O. Box 173, 1 Shizi Street\\
Suzhou, Jiangsu 215006\\
China\\
Email: jiangwenhua@suda.edu.cn

\bigskip
\noindent
 Cun-Hui Zhang \\
Department of Statistics and Biostatistics \\
Rutgers University\\
 Piscataway, NJ 08854 \\
 U.S.A.\\
 Email: cunhui@stat.rutgers.edu

\end{document}